\documentclass{siamltex}
\usepackage{amsmath,amsfonts,amssymb}
\usepackage{latexsym}
\usepackage{graphicx,overpic}
\usepackage{subfigure,caption}
\usepackage{pgfplots}
\pgfplotsset{compat=newest}
\usepackage{color}
\usepackage{booktabs}
\usepackage{caption}
\usepackage{lscape}
\usepackage{array}
\newcolumntype{C}[1]{>{\centering\let\newline\\\arraybackslash\hspace{0pt}}m{#1}}

\usepackage{tikz}

\newcommand{\R}{\mathbb R}
\newcommand{\bA}{\mathbf A}

\newcommand{\bH}{\mathbf H}
\newcommand{\bI}{\mathbf I}

\newcommand{\bP}{\mathbf P}

\newcommand{\bV}{\mathbf V}

\newcommand{\bn}{\mathbf n}
\newcommand{\be}{\mathbf e}

\newcommand{\bu}{\mathbf u}
\newcommand{\bU}{\mathbf U}
\newcommand{\bv}{\mathbf v}
\newcommand{\bw}{\mathbf w}

\newcommand{\bbf}{\mathbf f}

\newcommand{\T}{\mathcal T}

\newcommand{\divG}{{\mathop{\,\rm div}}_{\Gamma}}
\newcommand{\gradG}{\nabla_{\Gamma}}

\newcommand{\gradGh}{\nabla_{\Gamma_h}}

\newcommand{\OGamma}{\Omega^\Gamma_h}

\newcommand{\cT}{\mathcal T}

\newcommand{\vect}[1]{\boldsymbol{\mathbf{#1}}}

\newtheorem{remark}{Remark}[section]

\numberwithin{equation}{section}
\begin{document}

\title{Higher order Trace Finite Element Methods for the Surface Stokes Equation}
\author{Thomas Jankuhn\thanks{Institut f\"ur Geometrie und Praktische  Mathematik, RWTH-Aachen
University, D-52056 Aachen, Germany (jankuhn@igpm.rwth-aachen.de)} \and
Arnold Reusken\thanks{Institut f\"ur Geometrie und Praktische  Mathematik, RWTH-Aachen
University, D-52056 Aachen, Germany (reusken@igpm.rwth-aachen.de).}
}
\maketitle

\begin{abstract} 
In this paper a class of higher order finite element methods for the discretization of surface Stokes equations is studied. These methods are based on an unfitted finite element approach in which standard Taylor-Hood spaces on an underlying bulk mesh are used. For treating the constraint that the velocity must be tangential to the surface a penalty method is applied. Higher order geometry approximation is obtained by using a parametric trace finite element technique, known from the literature on trace finite element methods for scalar surface partial differential equations. Based on theoretical analyses for related problems, specific choices for the parameters in the method are proposed. Results of a systematic numerical study are included in which different variants are compared and convergence properties are illustrated. 
\end{abstract}
\section{Introduction}
In recent years there has been a strongly growing interest in the field of modeling and numerical simulation of surface fluids, cf. the papers \cite{arroyo2009,Jankuhn1,Kobaetal_QAM_2017,miura2017singular,Nitschkeetal_arXiv_2018,reuther2015interplay}, in which Navier-Stokes type PDEs on (evolving) surfaces  are treated. 
Navier-Stokes equations posed on manifolds is  a classical topic in  analysis, cf., e.g., \cite{ebin1970groups,mitrea2001navier,taylor1992analysis,Temam88}. There are only  very few papers that study numerical methods for surface (Navier-)Stokes equations \cite{nitschke2012finite,reuther2015interplay,reusken2018stream,reuther2018solving,fries2018higher,olshanskii2018finite,olshanskii2019penalty,Bonito2019a,OlshanskiiZhiliakov2019,Lederer2019}. Most of these papers consider either a (Navier-)Stokes system in stream function formulation (which assumes that the surface  is simlpy connected) or use low order $\vect P_1$-$P_1$ finite elements, combined with stabilization techniques. As far as we know, \cite{fries2018higher,OlshanskiiZhiliakov2019,Lederer2019} are the only papers in which  higher order finite element methods for surface Navier-Stokes equations are studied. In \cite{fries2018higher} a surface finite element approach \cite{DEreview} is used and the condition that the velocity must be tangential to the surface is enforced 
weakly by a Lagrange multiplier approach. In \cite{Lederer2019} a surface finite element approach is combined with a Piola transformation for the construction of divergence-free tangential finite elements. In \cite{OlshanskiiZhiliakov2019} stability of a variant of  the $\vect P_2$--$P_1$ Taylor-Hood pair for surface Stokes equations is analyzed and optimal discretization error bounds are derived. In this paper we consider a higher order finite element discretization that is based on a trace approach as in \cite{OlshanskiiZhiliakov2019}. For treating the tangential condition we use a penalty approach instead of the Lagrange multiplier method that is used in \cite{fries2018higher}. Instead of the surface finite element method of \cite{fries2018higher} we use a so-called trace finite element method. For scalar elliptic surface partial differential equations the surface and trace approaches are explained and compared in \cite{Bonito2019}. The former technique essentially uses an explicit surface triangulation with surface finite element spaces defined on it, whereas the latter approach uses an implicit (e.g., level set) representation of the surface combined with finite 
element spaces that are defined on an underlying unfitted bulk mesh.

The trace finite element method is a geometrically  \emph{un}fitted discretization. In particular in  a setting with
  evolving surfaces $\Gamma(t)$  such unfitted finite element techniques, also called cut FEM, have certain attractive properties concerning flexibility (no remeshing) and robustness (w.r.t.   handling of topological singularities); see \cite{olshanskii2016trace,cutFEM} for an overview of these techniques.

In this paper we study a  trace variant of the Taylor-Hood pair $\vect P_k$--$P_{k-1}$, $k \geq 2$, for discretization of surface Stokes equations. The case $k=2$ is treated in \cite{OlshanskiiZhiliakov2019}. Compared to Stokes equations in Euclidean domains, the surface variant leads to  several additional issues that have to be addressed. The two most important issues are the following:\\
1. \emph{Tangential flow constraint.} In surface flow problems the flow has to be tangential to the surface. It is not obvious how this  constraint (which is trivially satisfied in Euclidean domains) can be treated numerically. A technique used in several recent papers is as follows: the surface PDE for the tangential flow field is replaced by a PDE that allows fully three-dimensional velocities, defined on the surface, and a penalty approach is used to control the component of the velocity field that is normal to the surface.\\
2. \emph{Sufficiently accurate geometry approximation.} This topic resembles the problem of a sufficiently accurate boundary approximation for (Navier-)Stokes equations in Euclidean domains. For the latter the  isoparametric finite element technique is a standard approach. It is evident that for the case in which the domain is a curved surface the issue of geometry approximation becomes much more important.  To state it differently, for problems in Euclidean domains with a polygonal boundary, standard higher order finite elements (e.g., Taylor-Hood pair)  yield optimal higher order accuracy, whereas in a  finite element method for surface PDEs one \emph{always} needs a ``sufficiently accurate'' surface approximation for optimal higher order accuracy.

As mentioned above, we restrict to trace finite element techniques. Already for the case of scalar surface PDEs, in such trace methods  one applies an appropriate stabilization to control instabilities caused by ``small cuts''. 

In this setting of trace finite element methods for surface Stokes equations several important questions arise that are non-existent in Stokes problems in Euclidean domains. For example,  how does the error in the geometry approximation influence the discretization error, or,  what  is an appropriate scaling (in terms of the mesh size parameter $h$) of the penalty and stabilization parameters?

For the  trace variant of the Taylor-Hood pair that we present in this paper for all parameters, such as the order of polynomial degree used in the geometry approximation, the penalty parameter and stability parameters,  specific choices are proposed. These are based on analyses of related problems presented in \cite{jankuhn2019,OlshanskiiZhiliakov2019}. In \cite{jankuhn2019} an error analysis of a class of higher order trace finite element methods for a surface vector-Laplace problem is given. In \cite{OlshanskiiZhiliakov2019} the discrete inf-sup stability of a trace $\vect P_2$--$P_1$ Taylor-Hood pair is derived. 

Key ingredients of the higher order trace  finite element methods that we present in this paper are the following:
\begin{itemize}
 \item We use a penalty formulation for treating the tangential constraint. Two different variants will be studied,  namely a consistent and an inconsistent one. Precise explanations are given in Section~\ref{secttangential}.
 \item We use  parametric trace finite element spaces, known from scalar surface PDEs \cite{grande2017higher} and higher oder unfitted FEM for interface problems \cite{lehrenfeld2017analysis}, to obtain a higher order geometry approximation. The basic idea of this technique is outlined in Section~\ref{sectparametric}.
 \item The resulting trace finite element methods, including appropriate stabilization terms, are presented in Section~\ref{sectFEmethods}. The stabilization that we use, is the so-called \emph{volume normal derivative} stabilization, known from  the literature.
 \end{itemize}
 To decide on appropriate parameter choices, we briefly recall recently obtained rigorous stability and discretization error results for $\vect P_2$--$P_1$ surface Taylor-Hood elements  \cite{OlshanskiiZhiliakov2019} and  error bounds for trace FEM applied to surface vector-Laplace equations \cite{jankuhn2019}. These results are in given Section~\ref{sectionanalysis}. The proposed methods are applied to a surface Stokes equation on a sphere and on a torus. Results of numerical experiments that illustrate the optimal order of accuracy in different norms are presented in Section~\ref{sectionnumex}.

The topic of this paper relates to the one in \cite{OlshanskiiZhiliakov2019} as follows. In the latter paper the focus is on a theoretical analysis of discrete inf-sup stability of the trace $\vect P_2$--$P_1$ Taylor-Hood pair. An optimal order discretization error bound for this pair is derived in which, however, geometry errors are not treated. In this paper we focus on a general methodology for higher order trace $\vect P_k$--$P_{k-1}$ Taylor-Hood pairs, $k \geq 2$, in which the issue of geometry errors is also addressed. Furthermore, we compare two different penalty approaches, namely a consistent and an inconsistent one.

%We do not consider Lagrange multiplier technique (used in \cite{fries2018higher}).

\section{Continuous problem} \label{sectioncont}
We assume that $\Omega \subset \mathbb{R}^3$ is a polygonal domain which contains a connected compact smooth hypersurface $\Gamma$ without boundary. For the higher order finite element method that we introduce below it is essential that the surface $\Gamma$ is characterized as the zero level of a smooth level set function. For this we introduce some notation. A tubular neighborhood of $\Gamma$ is denoted by
$
U_\delta := \left\lbrace x \in \mathbb{R}^3 \mid \vert d(x) \vert < \delta \right\rbrace,
$
with $\delta > 0$ and $d$ the signed distance function to $\Gamma$, which we take negative in the interior of $\Gamma$. The surface $\Gamma$ is (implicitly) represented as the zero level of a smooth level set function $\phi \colon U_\delta \to \mathbb{R}$, i.e.
\[ \Gamma = \{ x \in \Omega \mid \phi(x) =0 \}.
\]
 This level set function is not necessarily close to a distance function but has the usual properties of a level set function:
\begin{equation*}
\Vert \nabla \phi(x) \Vert \sim 1, \quad \Vert \nabla^2 \phi(x) \Vert \leq c \quad \text{for all } x \in U_\delta.
\end{equation*}
We assume that the level set function   $\phi$  is sufficiently smooth. On $U_\delta$ we define $\bn(x) = \nabla d(x)$, the outward pointing unit normal on $\Gamma$, $\bH(x) = \nabla^2d(x)$, the Weingarten map,  $\bP = \bP(x):= \bI - \bn(x)\bn(x)^T$, the orthogonal projection onto the tangential plane, $p(x) = x - d(x)\bn(x)$, the closest point projection. We assume $\delta>0$ to be sufficiently small such that the decomposition $
x = p(x) + d(x) \bn(x)
$
is unique for all $x \in U_{\delta}$. The constant normal extension for vector functions $\bv \colon \Gamma \to \mathbb{R}^3$ is defined as $\bv^e(x) := \bv(p(x))$, $x \in U_{\delta}$. The extension for scalar functions is defined similarly. Note that on $\Gamma$ we have $\nabla \bv^e = \nabla(\bv \circ p) = \nabla \bv^e \bP$, with $\nabla \bw := (\nabla w_1, \nabla w_2, \nabla w_3)^T \in \mathbb{R}^{3 \times 3}$ for smooth vector functions $\bw \colon U_\delta \to \mathbb{R}^3$. For a scalar function $g \colon U_\delta \to \mathbb{R}$ and a vector function $\bv \colon U_\delta \to \mathbb{R}^3$ we define the surface (tangential and covariant) derivatives by
\begin{equation*} \begin{split}
\nabla_{\Gamma} g(x) &= \bP(x)\nabla g(x), \quad x \in \Gamma, \\
\nabla_{\Gamma} \bv(x) &= \bP(x)\nabla \bv(x) \bP(x), \quad x \in \Gamma.
\end{split}
\end{equation*} 
If $g$, $\bv$ are defined only on $\Gamma$, we use these definitions applied to the extension $g^e$, $\bv^e$.
%We also define the tangential derivative of a vector function $\bv \colon \Gamma \to \mathbb{R}^3$
%\begin{equation*}
%\nabla_{\Gamma,\tau} \bv(x) = \sum_{i=1}^3 \gradG \bv_i(x) = \nabla \bv^e(x) \bP(x), \quad x \in \Gamma.
%\end{equation*}
On $\Gamma$ we consider the surface stress tensor (see \cite{GurtinMurdoch75}) given by
\begin{equation*}
E_s(\bu):= \frac{1}{2} \left( \gradG \bu + \gradG^T \bu \right),
\end{equation*}
with $\gradG^T \bu := (\gradG \bu)^T$. To simplify the notation we write $E=E_s$.  The surface divergence operator for vector-valued functions $\bu \colon \Gamma \to \mathbb{R}^3$ and tensor-valued functions $\bA \colon \Gamma \to \mathbb{R}^{3\times3}$ are defined as
\begin{equation*} \begin{split}
\divG \bu &:= \textrm{tr} (\gradG \bu),  \\
\divG \bA &:= \left( \divG (\be_1^T\bA), \divG (\be_2^T\bA),\divG (\be_3^T\bA) ) \right)^T,
\end{split}
\end{equation*}
with $\be_i$ the $i$th basis vector in $\mathbb{R}^3$. For a given force vector $\bbf \in L^2(\Gamma)^3$, with $\bbf \cdot \bn =0$, and a source term $g \in L^2(\Gamma)$, with $\int_\Gamma g\, ds =0$, we consider the following \emph{surface Stokes problem}: determine $\bu \colon \Gamma \to \R^3$ with $\bu \cdot \bn = 0$ and $p \colon \Gamma \to \R$ with $\int_\Gamma p\, ds =0$ such that
\begin{equation} \label{eqstrong} \begin{split}
- \bP \divG (E(\bu)) + \bu + \gradG p &= \bbf \qquad \text{on } \Gamma,  \\
\divG \bu &= g \qquad \text{on } \Gamma.
\end{split}
\end{equation}
We added the zero order term on the left-hand side to avoid technical details related to the kernel of the strain tensor $E$ (the so-called Killing vector fields).
The surface Sobolev space of weakly differentiable vector valued functions is denoted by
\begin{equation} \label{eqdefH1}
\begin{gathered}
\bV:= H^1(\Gamma)^3, \quad \text{with} ~ \Vert \bu \Vert_{H^1(\Gamma)}^2 := \int_{\Gamma} \Vert \bu(s) \Vert_2^2 + \Vert \nabla \bu^e(s) \Vert_2^2 \, ds.
\end{gathered}
\end{equation}
%Note that $\|\nabla \bu^e\|_2=\|(\nabla \bu^e)^T\|_2$ and  on $\Gamma$ we have 
%\begin{equation} \label{defH1} (\nabla \bu^e)^T=\bP\big(\nabla u_1^e, \nabla u_2^e, \nabla u_3^e \big) = \sum_{i=1}^3 \gradG u_i \be_i^T .
%\end{equation}
% Hence, the norm in \eqref{eqdefH1} is a natural extension to vector valued functions of the usual scalar $H^1(\Gamma)$-norm. 
The corresponding subspace of \emph{tangential} vector field is denoted by
\begin{equation*}
\bV_T := \left\lbrace \bu \in \bV \mid \bu \cdot \bn =0 \right\rbrace.
\end{equation*}
A vector $\bu \in \bV$ can be orthogonally decomposed into a tangential and a normal part. We use the notation:
\begin{equation*}
\bu = \bP \bu + (\bu \cdot \bn)\bn = \bu_T + u_N\bn.
\end{equation*}
For $\textbf{u}, \textbf{v} \in \bV$ and $p \in L^2(\Gamma$) we introduce the bilinear forms
\begin{align}
a(\textbf{u}, \textbf{v}) &:= \int_\Gamma E(\bu) : E(\bv) \, ds + \int_\Gamma \textbf{u} \cdot \textbf{v} \, ds, \label{blfa} \\
b_T(\textbf{u}, p) &:= - \int_\Gamma p \divG \bu_T \, ds. \label{blfb} 
\end{align}
Note that in the definition of $b_T(\bu,p)$ only the \emph{tangential} component of $\bu$ is used, i.e., $b_T(\bu,p)=b_T(\bu_T,p)$ for all $\bu \in \bV$, $p\in L^2(\Gamma)$. This property motivates the notation $b_T(\cdot,\cdot)$ instead of $b(\cdot,\cdot)$.
If $p$ is from $H^1(\Gamma)$, then integration by parts yields
\begin{equation}\label{Bform}
b_T(\bu,p)=\int_\Gamma  \bu_T\cdot \gradG p \, ds = \int_\Gamma  \bu \cdot \gradG p \, ds.
\end{equation}
%For a given $\bbf \in L^2(\Gamma)^3$ with $\bbf \cdot \bn = 0$ and $g \in L^2(\Gamma)$ 
We introduce the following variational formulation of \eqref{eqstrong}: determine $(\bu_T, p)  \in \bV_T \times L^2_0(\Gamma)$ such that 
\begin{equation} \label{contform} \begin{split}
a(\bu_T, \bv_T) + b_T(\bv_T, p) &= (\bbf, \bv_T)_{L^2(\Gamma)} ~~~ \text{for all}~ \bv_T \in \bV_T, \\
b_T(\bu_T,q) &= (-g,q)_{L^2(\Gamma)} ~~~ \text{for all}~ q \in L^2(\Gamma).
\end{split}
\end{equation}
The bilinear form $a(\cdot,\cdot)$ is continuous on $\bV$, hence on $\bV_T$. The ellipticity of $a(\cdot,\cdot)$ on $\bV_T$ follows from the following surface Korn inequality, that holds if $\Gamma$ is $C^2$ smooth ((4.8) in \cite{Jankuhn1}): 
 There exists a constant $c_K > 0$ such that 
\begin{equation} \label{Korn}
\Vert \bu \Vert_{L^2(\Gamma)} + \Vert E(\bu) \Vert_{L^2(\Gamma)} \geq c_K \Vert \bu \Vert_{H^1(\Gamma)} \qquad \text{for all } \bu \in \bV_T.
\end{equation}
The bilinear form $b_T(\cdot,\cdot)$ is continuous on $\bV_T \times L_0^2(\Gamma)$ and satisfies the following inf-sup condition (Lemma 4.2 in \cite{Jankuhn1}):
There exists a constant $c>0$ such that estimate
\begin{equation*}
\inf_{p \in L^2_0(\Gamma)} \sup_{\bv_T \in \bV_T} \frac{b_T(\bv_T, p)}{\Vert \bv_T \Vert_{H^1(\Gamma)} \Vert p \Vert_{L^2(\Gamma)}} \geq c,
\end{equation*}
holds.
Hence, the weak formulation \eqref{contform} is \emph{a well-posed problem}. The unique solution is denoted by $(\bu_T^*,p^*)$. The main topic of this paper will be a class of higher order finite element methods for the discretization of this surface Stokes problem.

\section{Treatment of tangential flow constraint} \label{secttangential}
The weak formulation \eqref{contform} is not very suitable for a Galerkin finite element discretization, because we would need finite element functions that are (approximately)  tangential to $\Gamma$. Recently, very useful penalty approaches have been introduced \cite{Jankuhn1,hansbo2016analysis,jankuhn2019}. These techniques allow a full three-dimensional velocity $\bu$ (not necessarily tangential to $\Gamma$) and penalize the normal component of $\bu$. An alternative approach that avoids penalization has recently been introduced in \cite{Bonito2019a}.

In this section we recall two known penalty formulations: a consistent formulation and an inconsistent one. These are formulated as well-posed variational problems in a Hilbert space that contains $\bV_T$. In section~\ref{sectFEmethods} we apply a Galerkin discretization (modulo geometric errors) to these variational formulations. Both resulting finite element methods have there own merits, cf. Section~\ref{sectionconclusion}. 

We define $\bV_* \supset \bV \supset \bV_T$:
\begin{gather*}
\bV_* := \left\lbrace \bu \in L^2(\Gamma)^3 \mid \bu_T \in \bV_T, u_N \in L^2(\Gamma) \right\rbrace, 
\quad \Vert  \bu \Vert_{V_*}^2 := \Vert \bu_T \Vert_{H^1(\Gamma)}^2 + \Vert u_N \Vert_{L^2(\Gamma)}^2.
\end{gather*} 
Based on the  identity 
\begin{equation} \label{identity}
E(\bu) = E(\bu_T) + u_N \bH, \quad \bu \in \bV,
\end{equation}
we introduce an extension of the  bilinear form $a(\cdot,\cdot)$ from $\bV$ to the larger space $\bV_\ast$: 
\begin{equation} \label{defbla}
a(\textbf{u}, \textbf{v}) := \int_\Gamma (E(\bu_T) + u_N \bH) : (E(\bv_T) + v_N \bH) \, ds + \int_\Gamma \textbf{u} \cdot \textbf{v} \, ds, \quad \bu,\bv \in \bV_*.
\end{equation}
This bilinear form is well-defined and continuous on $\bV_*$. We also define a penalty bilinear form
\[
 k(\bu,\bv) := \eta \int_\Gamma  (\textbf{u} \cdot \textbf{n}) ~ (\textbf{v} \cdot \textbf{n})  \, ds \qquad \bu,\bv \in \bV_*,
\]
with $\eta >0$ a penalty parameter, and
\begin{equation*}
A(\bu,\bv):=a(\bu, \bv) + k(\bu, \bv) \qquad \bu,\bv \in \bV_*.
\end{equation*}We further introduce the bilinear form $a_T(\cdot,\cdot)$ in which only the tangential components of the arguments play a role:
\begin{equation} \label{bilinearformpenalty} 
a_T(\bu, \bv) := a(\bP\bu,\bP\bv)=a(\bu_T,\bv_T),
\end{equation}
and correspondingly,
\begin{equation*}
A_T(\bu,\bv):=a_T(\bu, \bv) + k(\bu, \bv) \qquad \bu,\bv \in \bV_*.
\end{equation*}
(Note that $A(\cdot,\cdot)$ and  $A_T(\cdot,\cdot)$ depend on the penalty parameter $\eta$).
A \emph{consistent} penalty surface Stokes formulation is: Determine $(\bu, p)  \in \bV_* \times L^2_0(\Gamma)$ such that 
\begin{equation} \tag{P1}\label{projectedcontform1} \begin{split}
A_T(\bu,\bv) + b_T(\bv, p) &= (\bbf, \bv)_{L^2(\Gamma)} ~~~ \text{for all}~ \bv \in \bV_*, \\
b_T(\bu,q) &= -(g,q)_{L^2(\Gamma)} ~~~ \text{for all}~ q \in L^2(\Gamma).
\end{split}
\end{equation} 
Using the surface Korn inequality \eqref{Korn} one obtains ellipticity of the bilinear form $A_T(\cdot, \cdot)$, which is used to derive  the following result (Theorem~6.1  in \cite{Jankuhn1}):
\begin{lemma}
Problem \eqref{projectedcontform1} is well-posed. For the unique solution $(\tilde{\bu} , \tilde{p}) \in \bV_* \times L^2_0(\Gamma)$ of this problem we have $(\tilde{\bu}, \tilde{p}) = (\bu_T^*, p^*)$.
\end{lemma}

The property
$(\tilde{\bu}, \tilde{p}) = (\bu_T^*, p^*)$ explains, why we call this a \emph{consistent} penalty formulation.

An \emph{inconsistent} penalty surface Stokes formulation is: Determine $(\bu, p)  \in \bV_* \times L^2_0(\Gamma)$ such that 
\begin{equation} \tag{P2}\label{projectedcontform2} \begin{split}
A(\bu,\bv) + b_T(\bv, p) &= (\bbf, \bv)_{L^2(\Gamma)} ~~~ \text{for all}~ \bv \in \bV_*, \\
b_T(\bu,q) &= -(g,q)_{L^2(\Gamma)} ~~~ \text{for all}~ q \in L^2(\Gamma).
\end{split}
\end{equation} 
From Theorem~3.1 in \cite{olshanskii2018finite} we get the following result:
\begin{lemma}
Assume $\eta$ is sufficiently large. Then  the problem \eqref{projectedcontform2} is well-posed and  for the unique solution $(\hat \bu, \hat p) \in \bV_* \times L^2_0(\Gamma)$ we have:
\begin{equation*}
\Vert \hat{\bu}_T - \bu_T^* \Vert_{H^1(\Gamma)} + \Vert \hat{u}_N \Vert_{L^2(\Gamma)} + \Vert \hat{p} - p^* \Vert_{L^2(\Gamma)} \leq C \eta^{-1}(\Vert \bbf \Vert_{L^2(\Gamma)} + \Vert g \Vert_{L^2(\Gamma)}).
\end{equation*}
\end{lemma}
The unique velocity solution $\hat \bu$ of \eqref{projectedcontform2} has a normal component that in general is nonzero. Due to $ \Vert \hat{u}_N \Vert_{L^2(\Gamma)} \leq C \eta^{-1}(\Vert \bbf \Vert_{L^2(\Gamma)} + \Vert g \Vert_{L^2(\Gamma)})$ its size can be controlled by the penalty parameter $\eta$.

\section{Parametric finite element space for high order surface approximation} \label{sectparametric}
%\textcolor{red}{$H^1(\Gamma_h)$ und Norm davon definieren.}
Clearly, for a higher order accurate finite element discretization  of the variational problems \eqref{projectedcontform1} and \eqref{projectedcontform2} one needs a sufficiently accurate approximation of the surface $\Gamma$. For this we use the \emph{parametric} trace finite element approach as in \cite{grande2017higher, jankuhn2019}. In this section we outline the  parametric mapping and the corresponding finite element space used in this method and summarize certain properties, known from the literature.

Let $\{ \mathcal{T}_h \}_{h>0}$ be a family of shape regular tetrahedral triangulations of $\Omega$. 
%For simplicity, in the analysis of the method, we assume $\{\T_h\}_{h >0}$ to be quasi-uniform. 
By $V_h^k$ we denote the standard finite element space of continuous piecewise polynomials of degree $k$. The nodal interpolation operator in $V_h^k$ is denoted by $I^k$. As input for the parametric mapping we need an approximation of $\phi$. We consider geometry approximations whose order of approximation may differ from the order of the polynomials used in the finite element space (introduced below). In other words, the spaces that we consider are not necessarily \emph{iso}parametric.  Let $k_g$ be the geometry approximation order, i.e., the construction of the geometry approximation will be based on a level set function approximation  $\phi_h \in V_h^{k_g}$. We assume that for this approximation the error estimate
\begin{equation} \label{phibound}
\max_{T \in \mathcal{T}_h} \vert \phi_h - \phi \vert_{W^{l,\infty}(T\cap U_\delta)} \leq c h^{k_g+1-l}, \quad 0\leq l \leq k_g+1,  
\end{equation}
is satisfied.
Here, $\vert \cdot \vert_{W^{l,\infty}(T\cap U_\delta)}$ denotes the usual semi-norm on the Sobolev space $W^{l,\infty}(T\cap U_\delta)$ and the constant $c$ depends on $\phi$ but is independent of $h$. The zero level set of the finite element function $\phi_h$ \emph{implicitly} characterizes an approximation of the interface, which, however, is hard to compute for $k_g \geq 2$.  With the piecewise \emph{linear} nodal interpolation of $\phi_h$, which is denoted by $\hat{\phi}_h = I^1\phi_h$, we define the low order geometry approximation:
\begin{equation*}
\Gamma^{\text{lin}} := \{ x \in \Omega \mid \hat{\phi}_h (x) = 0\},
\end{equation*}
which can easily be determined. 
The tetrahedra $T \in \mathcal{T}_h$ that have a nonzero intersection with $\Gamma^{\text{lin}}$ are collected in the set denoted by $\T_h^\Gamma$. The domain formed by all tetrahedra in $\T_h^\Gamma$ is denoted by $\OGamma:= \{ x \in T \mid T \in \T_h^\Gamma \}$. Let $\Theta_h^{k_g}\in  \big({V_h^{k_g}}_{|\OGamma}\big)^3$ be the mesh transformation of order $k_g$ as defined in \cite{grande2017higher}, cf. Remark~\ref{transfo}.
\begin{remark} \label{transfo} \rm We outline the key idea of the mesh transformation $\Theta_h^{k_g}$. For a detailed description and analysis we refer to \cite{grande2017higher,lehrenfeld2016high,lehrenfeld2017analysis}. There exists a unique $\tilde{d} \colon \Omega_h^\Gamma \to \mathbb{R}$ such that $\tilde{d}(x)$ is the in absolute value smallest number such that
\begin{equation*}
\phi\big(x + \tilde{d}(x) \nabla \phi(x)\big) = \hat{\phi}_h(x) \qquad \text{for } x \in \Omega_h^\Gamma.
\end{equation*}
Using $\tilde{d}$ we define the injective mapping 
\begin{equation*}
\Psi(x) := x + \tilde{d}(x) \nabla \phi(x), \qquad x \in \Omega_h^\Gamma,
\end{equation*}
which has  the property $\Psi(\Gamma^{\text{lin}}) = \Gamma$. This mapping $\Psi$ deforms the mesh $\T_h^\Gamma$ in such a way that the (available) surface approximation $\Gamma^{\text{lin}}$ is mapped to the exact surface $\Gamma$.  To avoid computations with $\phi$ (which even may not be available) we use a similar construction with $\phi$ replaced by its (finite element) approximation $\phi_h$.
The resulting mapping $\Psi_h$ is not necessarily a finite element function. The  mesh transformation $\Theta_h^{k_g}$ is obtained by a simple projection (based on local averaging of values around a vertex) of $\Psi_h$ into the finite element space $\big({V_h^{k_g}}_{|\OGamma}\big)^3$. This parametric mapping is easy to determine. Implementation aspects are discussed in \cite{lehrenfeld2016high}. The mapping is implemented in Netgen/NGSolve \cite{ngsolve}.
\end{remark}

An approximation of $\Gamma$ is defined by
\begin{equation*}
\Gamma_{h}^{k_g} := \Theta_h^{k_g}(\Gamma^{\text{lin}}) = \left\lbrace x \mid \hat{\phi}_h((\Theta_{h}^{k_g})^{-1}(x)) = 0 \right\rbrace.
\end{equation*}
In \cite{lehrenfeld2017analysis} it is shown that (under certain reasonable smoothness assumptions) the estimate
\begin{equation} \label{distres}
  {\rm dist}(\Gamma_{h}^{k_g}, \Gamma) \lesssim h^{k_g+1}
\end{equation}
holds. Here and further in the paper we write $x \lesssim y$ to state that there exists a constant $c>0$, which is independent of the mesh parameter $h$ and the position of $\Gamma$ in the background mesh, such that the inequality $x \leq cy$ holds. Hence, the paramatric mapping $\Theta_{h}^{k_g}$ indeed yields a higher order surface approximation. 
We denote the transformed cut mesh domain by $\Omega^{\Gamma}_\Theta := \Theta_h^{k_g}(\Omega^\Gamma_h)$ and apply to $V_h^{k}$ the transformation $\Theta_h^{k_g}$ resulting in the \emph{parametric spaces} (defined on $\Omega^{\Gamma}_\Theta$)
\begin{equation*} 
V_{h,\Theta}^{k,k_g} := \left\lbrace v_h \circ (\Theta_h^{k_g})^{-1} \mid v_h \in {V_h^{k}}_{|\Omega^{\Gamma}_h} \right\rbrace, \quad
\bV_{h,\Theta}^{k,k_g} := (V_{h,\Theta}^{k,k_g})^3.
\end{equation*}
Note that $k_g$ denotes the degree of the polynomials used in the parametric mapping $\Theta_h^{k_g}$, which determines the accuracy of the  geometry approximation, cf.~\eqref{distres}, and $k$ the degree of the polynomials used in the finite element space. 
To simplify the notation we delete the superscript $k_g$ and write 
\[ V_{h,\Theta}^{k} = V_{h,\Theta}^{k,k_g},~~\bV_{h,\Theta}^{k} = \bV_{h,\Theta}^{k,k_g},~~\Theta_h = \Theta_h^{k_g}, ~~\Gamma_{h} = \Gamma_{h}^{k_g}.\]  

The following lemma, taken from \cite{grande2017higher}, gives an approximation error for the easy to compute normal approximation $\bn_h$, which is used in the methods introduced below.
\begin{lemma} \label{lemmanormals}
For $x \in T \in \mathcal{T}^\Gamma_h$ define
\begin{equation*}
\bn_{\textrm{lin}} = \bn_{\textrm{lin}}(T) := \frac{\nabla \hat{\phi}_h(x)}{\Vert \nabla \hat{\phi}_h(x)\Vert_2} = \frac{\nabla \hat{\phi}_{h|T}}{\Vert \nabla \hat{\phi}_{h|T}\Vert_2}, \quad \bn_h (\Theta(x)) := \frac{D\Theta_h(x)^{-T} \bn_{\textrm{lin}}}{\Vert D\Theta_h(x)^{-T} \bn_{\textrm{lin}} \Vert_2}.
\end{equation*}
Let $\bn_{\Gamma_h}(x)$, $x \in \Gamma_h$ a.e., be the unit normal on $\Gamma_h$ (in the direction of $\phi_h >0$). The following holds:
\begin{equation*} \begin{split}
\Vert \bn_h - \bn \Vert_{L^\infty(\Omega_\Theta^\Gamma)} &\lesssim h^{k_g}, \\
\Vert \bn_{\Gamma_h} - \bn \Vert_{L^\infty(\Gamma_h)} &\lesssim h^{k_g}.
\end{split}
\end{equation*}
\end{lemma}
%Similarly for $x\gtrsim y$, and $x \sim y$  means that both $x\lesssim y$ and $x\gtrsim y$ hold.
%Similar to the extension of a function $u$ defined on $\Gamma$ to $u^e$ defined on $U_\delta$ we define the lifting $u^l$ of a function $u$ defined   on $\Gamma_h$ by 
%\begin{equation*}
%\begin{cases}
%u^l(p(x)) = u(x) &\text{ for } x \in \Gamma_h, \\
%u^l(x) = u^l(p(x)) &\text{ for } x \in U_\delta.
%\end{cases}
%\end{equation*} 
%A norm on $H^1(\Gamma_h)^3$ is defined using the component-wise lifting by
%\begin{equation*}
%\Vert \bu \Vert_{H^1(\Gamma_h)}^2 := \int_{\Gamma_h} \Vert \bu(s) \Vert_2^2 + \Vert \nabla \bu^l(s) \bP_h(s) \Vert_2^2 \, ds
%\end{equation*}
%with $\bP_h = \bI - \bn_h \bn_h^T$. 
%%Similar to the definition of the norm on $H^1(\Gamma)^3$ we only consider gradients of the components of $\bu$ which are tangential to $\Gamma_h$. 
%In \eqref{eqdefH1} the constant normal extension automatically yields tangential gradients of all components, cf. \eqref{defH1}.
% The lifting used in the definition of the $H^1(\Gamma_h)$-norm is constant along the normal to $\Gamma$ (not $\Gamma_h$). Therefore,  to eliminate the part of the (componentwise) gradient which is normal to $\Gamma_h$ one needs the projection $\bP_h$. 
%We also introduce the  following  spaces
%\begin{equation*} \begin{split}
%V_{reg,h} := \left\lbrace v \in H^1(\Omega_{\Theta}^{\Gamma}) \mid \text{tr}_{|_{\Gamma_h}} v \in H^1(\Gamma_h) \right\rbrace \supset V_{h,\Theta}^{k}, \quad
%\bV_{reg,h} := \big(V_{reg,h}\big)^3.
%\end{split}
%\end{equation*}
%
\section{Higher order  trace finite element methods} \label{sectFEmethods}
In this section we introduce a class of higher order parametric trace finite element methods. These methods  are obtained by applying a Galerkin approach (modulo a geometry error due to $\Gamma_h \approx \Gamma$) to the formulations \eqref{projectedcontform1} and \eqref{projectedcontform2}. Based on the parametric finite element spaces $\bV_{h,\Theta}^{k}$ and $V_{h,\Theta}^{k}$ we introduce for $k \geq 2$ the $\vect P_k$-$P_{k-1}$ pair of \emph{parametric trace Taylor-Hood elements}:
\begin{equation*}
\bU_h := \bV_{h,\Theta}^{k}, \qquad Q_h := V_{h,\Theta}^{k-1} \cap L^2_0(\Gamma_h).
\end{equation*}
Note that  the polynomial degrees, $k$ and $k-1$, for the velocity and pressure approximation are different, but both spaces $\bU_h$ and $Q_h$ use the same parametric mapping based on polynomials of degree $k_g$. Since the pressure approximation uses $H^1$ finite element functions we can use the partial integration \eqref{Bform} (with $\Gamma$ replaced by $\Gamma_h$).
We introduce    discrete variants of the bilinear forms $a(\cdot,\cdot)$, $a_T(\cdot,\cdot)$, $b_T(\cdot,\cdot)$ and the penalty bilinear form $k(\cdot,\cdot)$ introduced above.  Since we use a trace FEM, we need a stabilization that eliminates instabilities caused by the small cuts. For this we use the so-called ``normal derivative volume stabilization'', known from the literature \cite{burmanembedded,grande2017higher}  ($s_h(\cdot,\cdot)$ and $\tilde{s}_h(\cdot,\cdot)$ below). We define, with $\bP_h=\bP_h(x):=\bI- \bn_h(x)\bn_h(x)^T$, $x \in \Omega_{\Theta}^{\Gamma}$:
\begin{align*}
 \gradGh \bu &:= \bP_h \nabla \bu \bP_h, \\  E_h(\bu) &:= \frac12 \big(\gradGh \bu + \gradGh^T \bu\big), \quad E_{T,h}(\bu):=E_h(\bu) - u_N \bH_h, \\
a_{h}(\bu,\bv) &:= \int_{\Gamma_h} E_{h}(\bu):E_{h}(\bv)\, ds_h + \int_{\Gamma_h} \bu \cdot \bv \, ds_h,\\
 a_{T,h}(\bu,\bv) &:= \int_{\Gamma_h} E_{T,h}(\bu):E_{T,h}(\bv)\, ds_h + \int_{\Gamma_h} \bP_h \bu \cdot \bP_h \bv \, ds_h,\\ 
  b_h(\bu,q)& := \int_{\Gamma_h} \bu \cdot \gradGh q \, ds_h,\\
 k_h(\bu,\bv)&:= \eta \int_{\Gamma_h} (\bu \cdot \tilde{\bn}_h) (\bv \cdot \tilde{\bn}_h)  \, ds_h, \\
  s_h(\bu,\bv) & := \rho_u\int_{\Omega_{\Theta}^{\Gamma}} (\nabla \bu \bn_h) \cdot (\nabla \bv \bn_h)  \, dx, \quad
  \tilde{s}_h(p,q) := \rho_p  \int_{\Omega_{\Theta}^{\Gamma}} (\bn_h\cdot \nabla p) (\bn_h \cdot \nabla q) \, dx.
\end{align*}
%All these bilinear forms are well-defined for $\bu,\bv \in \bV_{reg,h}$ and $p,q  \in V_{reg,h}$.
The normal vector $\tilde{\bn}_h$, used in the penalty term $k_h(\cdot, \cdot)$, and the curvature tensor $\bH_h$ are approximations of the exact normal and the exact Weingarten mapping, respectively.
The reason that we introduce yet another normal approximation $\tilde{\bn}_h$ is the following. From an error  analysis of the vector-Laplace problem in \cite{hansbo2016analysis,jankuhn2019}, cf. also section~\ref{sectionanalysis} below, it follows that for obtaining optimal order estimates  the normal approximation $\tilde{\bn}_h$ used in the penalty term has to be more accurate than the normal approximation $\bn_h$. How suitable approximations $\tilde{\bn}_h$ and $\bH_h$ can be determined is discussed in Section~\ref{sectionanalysis}.
Suitable choices of the stabilization parameters $\rho_u$, $\rho_p$ and the penalty parameter $\eta$ are  also discussed in  Section~\ref{sectionanalysis}.  
As a discrete analogon of  $E(\bu_T)= E(\bP\bu)= E(\bu)-u_N \bH$ we  use $E_{T,h}(\bu)=E_h(\bu) - u_N \bH_h$ instead of $E_{T,h}(\bu)=E_h(\bP_h\bu)$, because the latter requires (tangential) differentiation of $\bP_h$, which causes difficulties. The (canonical) choice of $\bn_h$ as in Lemma~\ref{lemmanormals} is discontinuous across faces, hence \emph{not} an $H^1(\Gamma_h)$ vector function, which implies that $E_h(\bP_h \bu)$ is in general not well-defined. 
%The reason that we introduce yet another normal approximation $\tilde{\bn}_h$ is the following. From an error  analysis of the vector-Laplace problem in \cite{jankuhn2019}, cf. also section~\ref{sectionanalysis} below, it follows that for obtaining optimal order estimates  the normal approximation $\tilde{\bn}_h$ used in the penalty term has to be more accurate than the normal approximation $\bn_h$. 
%To get a consistent continuous formulation we need the additional tangential projections $\bP$ in the bilinear form of the differential operator, which results as we have seen above in the approximation of the Weingarten map in the discrete formulation. However, omitting that approximation of the Weingarten map in the discrete formulation leads to a formulation which does not have optimal high order discretization error bounds (see \cite{jankuhn2019}). 
%Because we take the trace of outer finite element functions on the surface approximation $\Gamma_h$, there is the issue that there may be two different finite element functions with the same trace on $\Gamma_h$. To deal with this non-uniqueness issue we stabilize the variational formulation using the "normal derivative volume" stabilization $s_h(\cdot,\cdot)$, which was introduced in \cite{grande2017higher}. 
%The stabilization with  $s_h(\cdot,\cdot)$ used in the variational penalty formulations below guarantees that the stiffness matrix has a spectral condition number $\sim h^{-2}$, independent of how the interface cuts the outer triangulation.
%\end{remark}

We now introduce  discrete versions of the formulations \eqref{projectedcontform1}  and \eqref{projectedcontform2}. For these we  need a suitable (sufficiently accurate) extension of the data $\bbf$ and $g$ to $\Gamma_h$, which are denoted by $\bbf_h$ and $g_h$, respectively. 
\\[1ex]
\emph{Consistent discrete surface Stokes.}  This method is based on the \emph{consistent} formulation \eqref{projectedcontform1} and uses the bilinear form $a_{T,h}(\cdot,\cdot)$. Define
\[
  A_{T,h}(\bu,\bv) := a_{T,h}(\bu,\bv) + s_h(\bu,\bv)+k_h(\bu,\bv).
\]
The discrete problem reads:
 determine $(\bu_h, p_h) \in \bU_h \times Q_h$ such that
\begin{equation} \tag{P1h} \label{discreteform1}
\begin{aligned}
A_{T,h}(\bu_h,\bv_h) + b_h(\bv_h,p_h) & =(\mathbf{f}_h,\bv_h)_{L^2(\Gamma_h)} &\quad &\text{for all } \bv_h \in \bU_h \\ 
b_h(\bu_h,q_h) - \tilde{s}_h(p_h,q_h) & = (-g_h,q_h)_{L^2(\Gamma_h)} &\quad &\text{for all }q_h \in Q_h.
\end{aligned}
\end{equation}
Note that, although we call this method ``consistent'', due to geometry errors it does contain  consistency errors.\\[1ex]
\emph{Inconsistent discrete surface Stokes.}
 This method is based on the \emph{inconsistent} formulation \eqref{projectedcontform2} and uses the bilinear form $a_h(\cdot,\cdot)$.
Define
\[ A_{h}(\bu,\bv) := a_{h}(\bu,\bv) + s_h(\bu,\bv)+k_h(\bu,\bv). \]
The discrete problem reads:
 determine $(\bu_h, p_h) \in \bU_h \times Q_h$ such that
\begin{equation} \tag{P2h} \label{discreteform2}
\begin{aligned}
A_h(\bu_h,\bv_h) + b_h(\bv_h,p_h) & =(\mathbf{f}_h,\bv_h)_{L^2(\Gamma_h)} &\quad &\text{for all } \bv_h \in \bU_h \\ 
b_h(\bu_h,q_h) - \tilde{s}_h(p_h,q_h) & = (-g_h,q_h)_{L^2(\Gamma_h)} &\quad &\text{for all }q_h \in Q_h.
\end{aligned}
\end{equation}

In the next section we explain how components of these methods, for example the penalty parameter $\eta$ and the Weingarten mapping approximaiton $\bH_h$, can be chosen. In Section~\ref{sectionnumex} we present numerical results for these methods.

\section{Choice of method components based on available analysis} \label{sectionanalysis}
Before the  finite element discretizations~\eqref{discreteform1} and \eqref{discreteform2} can be applied to a specific surface Stokes problem, the following issues have to be addressed:
\begin{itemize}
 \item[a)]  Accuracy of geometry approximation: given $k$, how should one take $k_g$?
 \item[b)] Components in penalty term: how does $\eta$ depend on $h$, how to choose $\tilde \bn_h$?
 \item[c)] Parameters in volume normal derivative stabilizations: how do $\rho_u$, $\rho_p$ depend on $h$?
 \item[d)] Weingarten mapping approximation (only for consistent method): what is a suitable choice for $\bH_h$?
\end{itemize}
In this section we address these issues and give specific recommendations. These are based on recent analyses of surface vector-Laplace and surface Stokes equations. Below we first summarize a few relevant results of these analyses that will be used to answer the questions above. It is convenient to introduce one  further order parameter $k_p \geq k$ (besides $k$ and $k_g$) that describes the accuracy of the normal approximation $\tilde \bn_h$:
\begin{equation} \label{apprbnt}
\Vert \bn - \tilde{\bn}_h \Vert_{L^{\infty}(\Gamma_h)} \lesssim h^{k_p}.
\end{equation}
In \cite{jankuhn2019} discrete vector-Laplace problems are studied that can be seen as  simplifications of the problems \eqref{discreteform1} and \eqref{discreteform2}. More precisely, in the vector-Laplace equation, the only unknown is a tangential velocity field $\bu$ (no pressure) that has to satisfy the  equation $- \bP \divG (E(\bu)) + \bu = \bbf$ on $\Gamma$, which is similar to \eqref{eqstrong}. The same parametric finite element techniques as described above are applied and yield discrete problems as in the first equations in \eqref{discreteform1} and \eqref{discreteform2}, with $b_h(\cdot,\cdot)$ put to zero. For these discretizations a complete error analysis (including geometry errors) is presented in \cite{jankuhn2019}. In that analysis the natural energy norm $\|\cdot \|_A$, defined by $\|\bv\|_A^2=\Vert \bv \Vert_{A_{T,h}}^2$ for the consistent method and $\|\bv\|_A^2=\Vert \bv \Vert_{A_{h}}^2$ for the inconsistent one, is used. Main results of the error analysis are the following (we refer to \
cite{jankuhn2019} for precise formulations of these results):
\begin{itemize}
 \item For the consistent method. Assume $\Vert \bH - \bH_h \Vert_{L^{\infty}(\Gamma_h)} \lesssim h^{k_{g}-1}$, $k_g=k$ (isoparametric case), $\eta \sim h^{-2}$, $k_p=k+1$, $\rho_u\sim h^{-1}$, $\rho_p \sim h$. Then an  optimal order error bound of order $\mathcal{O}(h^{k})$ in the energy norm holds. This bound implies an optimal error bound in the $H^1(\Gamma_h)$-norm of the same order.
\item For the inconsistent method.  Assume $k_g=k$ (isoparametric case), $\eta \sim h^{-(k+1)}$, $k_p=k+1$, $\rho_u\sim h^{-1}$, $\rho_p \sim h$. Then an  optimal order error bound of order $\mathcal{O}(h^{\frac12(k+1)})$ in the energy norm holds. This bound implies an error bound in the $H^1(\Gamma_h)$-norm of the same order, which is optimal only for the case $k=1$.
\end{itemize}
Furthermore, numerical experiments indicate the following:
\begin{itemize}
\item The inconsistent method, with parameters as above, has optimal order $\mathcal{O}(h^{k})$-convergence in the $H^1(\Gamma_h)$-norm and optimal order $\mathcal{O}(h^{k+1})$-convergence in the $L^2(\Gamma_h)$-norm  not only for $k=1$ but also for $k \geq 2$. 
\item Taking $k_p=k$ leads to suboptimal convergence behavior for both the consistent and the inconsistent method.
\item For the inconsistent method and $k \geq 2$, optimal order convergence is lost if for the penalty parameter we use a scaling $\eta \sim h^{-2}$.
%\item For optimal order of convergence of the  consistent method it is sufficient to have an approximation of the Weingarten mapping with accuracy $\Vert \bH - \bH_h \Vert_{L^{\infty}(\Gamma_h)} \lesssim h^{k_{g}-1}$. 
\end{itemize}
For these results to hold, one needs a sufficiently accurate data extension $\mathbf{f}_h$ of $\mathbf{f}$. Precise conditions are given in \cite{jankuhn2019} and are very similar to the conditions formulated for higher order methods for \emph{scalar} surface PDEs \cite{demlow2009higher,reusken2015analysis}.

In the recent paper \cite{OlshanskiiZhiliakov2019} the discretizations \eqref{discreteform1} and \eqref{discreteform2} are studied for the case \emph{without geometry errors}, i.e., $\Gamma_h=\Gamma$. In that case we do not need the paramatric mapping $\Theta_h$ and the finite element spaces are simply the Taylor-Hood pairs on the local triangulation, consisting of the tetrahedra intersected by $\Gamma$. Clearly, this method is in general not feasible in practice, because integrals over $\Gamma$ can  not be evaluated efficiently. This (simplified) setting, however, is used to analyze the discrete inf-sup stability of the trace Taylor-Hood pair for the surface Stokes problem. A main result derived in \cite{OlshanskiiZhiliakov2019} is the following (we refer to \cite{OlshanskiiZhiliakov2019} for precise formulation):
\begin{itemize}
 \item Assume $h \lesssim \rho_u \lesssim h^{-1}$, $\rho_p \sim h$, $\eta \sim h^{-2}$. Then both for the consistent and inconsistent variant the discrete inf-sup stability estimate
\begin{equation} \label{discrinf}  \|q\|_{L^2(\Gamma)}\lesssim \sup_{\bv\in\bU_h}\frac{b_T(\bv,q)}{\|\bv\|_A} + \tilde s_h(q,q)^{\frac12}\quad\text{for all}~q\in Q_h,
\end{equation}
holds for $k=2$, i.e., for the $\vect P_2$--$P_1$ trace Taylor-Hood pair. 
\item For this parameter choice of  $\rho_u$, $\rho_p$ and $\eta$ the  consistent method has an optimal error bound (in $H^1(\Gamma)$-norm for velocity and $L^2(\Gamma)$-norm for pressure). 
\end{itemize}
Based on these results, for the discretizations \eqref{discreteform1} and \eqref{discreteform2} of the surface Stokes problem we have the following recommendations concerning the issues a)-d) raised above.\\
a) \emph{ Accuracy of geometry approximation.} We take $k_g=k$, i.e. isoparametric finite elements for velocity.\\ 
b) \emph{ Components in penalty term.} For the consistent method we take $\eta \sim h^{-2}$ and for the inconsistent method $\eta \sim h^{-(k+1)}$. In both methods we  use a normal approximation $\tilde \bn_h$ with accuracy $k_p=k+1$. Such an approximation can be determined as follows. We assume that we have an approximation $\tilde \phi_h$ of $\phi$ available that is one order more accurate than $\phi_h$, i.e., it satisfies an error bound as in \eqref{phibound} with $k_g$ replaced by $k_g+1$. We then take $\tilde \bn_h:= \frac{\nabla \tilde \phi_h}{\|\nabla \tilde \phi_h\|}$. 
\\
c) \emph{Parameters in volume normal derivative stabilizations.} We take  $\rho_u \sim h^{-1}$, $\rho_p \sim h$. 
\\
d) \emph{Weingarten mapping approximation} (only for consistent method).  We use an approximation  $\bH_h$ with order of accuracy $k_g-1$. 
Such an approximation can be obtained  by taking $\bH_h = \nabla(I_{\Theta}^{k_g}(\bn_h))$, where $I_{\Theta}^{k_g}$ denotes the (componentwise) parametric nodal interpolation in the space $V_{h, \Theta}^{k_g}$, cf. \cite{grande2017higher}.

\section{Numerical experiments} \label{sectionnumex}

In this section we present results of numerical experiments. As test cases we consider Stokes equations on a sphere and a torus. For these two cases we first describe the setting of the continuous problem.

The unit sphere $\Gamma$ is characterized by the zero level of the distance function $\phi(x) = \sqrt{x_1^2+x_2^2+x_3^2} - 1$, $x= (x_1,x_2,x_3)^T$.
% i.e. $\Gamma = \{ x \in \Omega \mid \phi(x) =0 \}$. 
The surface is embedded in the domain $\Omega = [-5/3,5/3]^3$.  We consider the surface Stokes problem \eqref{contform} with the prescribed solution
\begin{equation*} \begin{split}
\bu(x) &= \begin{pmatrix} \frac{(x_2^2 x_3^2 + x_3^4) \sqrt{x_1^2+x_2^2+x_3^2} + x_1 (x_1^2+x_2^2+x_3^2) (x_1 x_3 + x_2^2)}{(x_1^2+x_2^2+x_3^2)^\frac52} \\ \frac{(x_1 x_3^2 \sqrt{x_1^2+x_2^2+x_3^2} +(x_1^2-x_1 x_3+x_3^2) (x_1^2+x_2^2+x_3^2))x_2}{(x_1^2+x_2^2+x_3^2)^\frac52} \\ \frac{x_1x_3^3 \sqrt{x_1^2+x_2^2+x_3^2}+(x_1^2+x_2^2+x_3^2)(x_1^3+x_1 x_2^2-x_2^2 x_3)}{(x_1^2+x_2^2+x_3^2)^\frac52} \end{pmatrix}, \\
p(x) &= \frac{x_1 x_2^3 + x_3 (x_1^2+x_2^2+x_3^2)^\frac{3}{2}}{(x_1^2+x_2^2+x_3^2)^2}.
\end{split}
\end{equation*}
The velocity solution is tangential, i.e. $\bP\bu = \bu$ and constant in normal direction, i.e., $\bu = \bu^e$. The velocity field $\bu$ is not divergence-free, i.e., $\divG \bu \neq 0$. The pressure solution is also constant in normal direction, i.e. $p = p^e$ as well as mean free, i.e. $\int_{\Gamma} p \, ds = 0$. Corresponding  right-hand sides $\bbf$ and $g$ are computed in a small neighborhood of $\Gamma$ as follows. The surface  differential operators used in the Stokes problem \eqref{eqstrong}, defined on $\Gamma$, have canonical extensions to a small neighborhood of $\Gamma$. We use these extended ones and  apply the Stokes operator (defined in the neighborhood) to the prescribed $\bu$ and $p$, which are constant in normal direction. The resulting $\bbf$ and $g$, which are defined in the neighborhood and not necessarily constant in normal direction,  are used as data $\bbf_h$ and $g_h$ in the finite element method. 
\\[1ex]
For the case of a torus, $\Gamma$ is characterized by the zero level of the distance function $\phi(x) = \sqrt{x_3^2 + (\sqrt{x_1^2+x_2^2} - 1)^2} - \frac{1}{2}$. The surface is again embedded in the domain $\Omega = [-5/3,5/3]^3$.  We consider the surface Stokes problem \eqref{contform} with the prescribed solution
\begin{equation*} \begin{split}
\bu(x) &= \bv^e(x) \qquad \text{with} \quad \bv(x) =  \begin{pmatrix} {\frac {x_{{3}}^{2}x_{{1}}}{ \left( x_{{1
}}^{2}+x_{{2}}^{2}+x_{{3}}^{2}-2\,\sqrt {x_{{1}}^{2}+x_{{2}}^
{2}}+1 \right) \sqrt {x_{{1}}^{2}+x_{{2}}^{2}}}}
\\ \noalign{\medskip}{\frac {x_{{2}}x_{{3}}^{2}}{ \left( x_{{1}}^{
2}+x_{{2}}^{2}+x_{{3}}^{2}-2\,\sqrt {x_{{1}}^{2}+x_{{2}}^{2}}+
1 \right) \sqrt {x_{{1}}^{2}+x_{{2}}^{2}}}}\\ \noalign{\medskip}-{
\frac { \left( \sqrt {x_{{1}}^{2}+x_{{2}}^{2}}-1 \right) x_{{3}}}{
x_{{1}}^{2}+x_{{2}}^{2}+x_{{3}}^{2}-2\,\sqrt {x_{{1}}^{2}+x_{
{2}}^{2}}+1}} \end{pmatrix}, \\
p(x) &= q^e(x) - \frac{\int_{\Gamma} q^e \, ds}{\int_{\Gamma} 1 \, ds} \qquad \text{with} \quad  q(x) = x_1 x_2^3+x_3.
\end{split}
\end{equation*}
The velocity solution is tangential, i.e. $\bP\bu = \bu$ and constant in normal direction, i.e. $\bu = \bu^e$. The velocity field $\bu$ is not divergence-free, i.e., $\divG \bu \neq 0$. The pressure solution is also constant in normal direction, i.e., $p = p^e$. The right-hand sides $\bbf_h$ and $g_h$ for the finite element discretization  are computed in the same way as for the sphere above. 

In both cases, for the construction of the local triangulation $\cT_h^\Gamma$
we start with an unstructured tetrahedral Netgen-mesh with $h_{max} = 0.5$ (see \cite{Schoeberl1997}) and locally refine the mesh  using a marked-edge bisection method (refinement of tetrahedra that are intersected by the surface).

In the implementation of the discretizations~\eqref{discreteform1} and \eqref{discreteform2} of the surface Stokes problem, we use (unless stated otherwise) the parameter setting and components listed in a)-d) at the end of section~\ref{sectionanalysis} (with a constant 1 in $\sim$). The methods are implemented in Netgen/NGSolve with ngsxfem \cite{ngsolve,ngsxfem}.

The errors are measured in different (semi-)norms. We use the following notations:
\begin{align*}
e_{L^2}^{\bu} &:= \Vert \bu - \bu_h \Vert_{L^2(\Gamma_h)}, & e_{H^1}^{\bu} &:= \Vert \gradGh(\bu - \bu_h) \Vert_{L^2(\Gamma_h)}, & e_{M}^{p} &:= \Vert p - p_h \Vert_{M}, \\ 
e_{PL^2}^{\bu} &:= \Vert \bP_h(\bu - \bu_h)\Vert_{L^2(\Gamma_h)}, & e_{A}^{\bu} &:= \Vert \bu - \bu_h \Vert_{A}.
\end{align*}
Here $\|\cdot\|_M^2:= \Vert \cdot \Vert_{L^2(\Gamma_h)}^2 + \tilde{s}_h(\cdot,\cdot)$.

\subsection{Results for the sphere}

In Section \ref{sectionoptimalresultssphere} we present results of numerical experiments that show optimal convergence orders for $\vect P_2$-$P_1$ and $\vect P_3$-$P_2$ finite elements and comment on the choice of the stabilization parameters $\rho_u$ and $\rho_p$. We also compare the consistent and inconsistent methods. In Section \ref{sectionpenaltyterm} we discuss the choice of the parameters in the penalty term and the effects the  penalty term has on the  energy norm error.

\subsubsection{Optimal results for $\vect P_2$-$P_1$ and $\vect P_3$-$P_2$ finite elements} \label{sectionoptimalresultssphere}
\begin{figure}
\hspace*{-0.7cm}
\begin{minipage}{0.51\textwidth}
\scriptsize
      % \begin{figure}  
\scalebox{0.81}{
  \begin{tikzpicture}
  \def\vara{0.18}
  \def\varb{0.25}
  \def\varc{0.1}
\begin{semilogyaxis}[ xlabel={Refinement level}, ylabel={Error}, ymin=5E-7, ymax=3, legend style={ at={(0.5,1.02)}, anchor=south, legend columns =2, transpose legend }, legend cell align=left, cycle list name=mark list, title= {$\vect P_2$-$P_{1}$}, title style={at={(0.5,1.16)}, align=center, font=\Large} ]
    \addplot table[x=level, y=L2] {Ph1P2P1opt.dat};
    \addplot table[x=level, y=L2P] {Ph1P2P1opt.dat};
    \addplot table[x=level, y=realH1semi] {Ph1P2P1opt.dat};
    \addplot table[x=level, y=uerror] {Ph1P2P1opt.dat};
    \addplot table[x=level, y=merr] {Ph1P2P1opt.dat};
%    \addplot[dashed,line width=0.75pt] coordinates { % h
%    (1,\vara) (2,\vara*0.5) (3,\vara*0.25) (4,\vara*0.125) (5,\vara*0.0625)(6,\vara*0.03125)
%    };
    \addplot[dashed,line width=0.75pt] coordinates { % h^2
      (1,\varb) (2,\varb*0.5*0.5) (3,\varb*0.25*0.25) (4,\varb*0.125*0.125) (5,\varb*0.0625*0.0625)(6,\varb*0.03125*0.03125)
    };
    \addplot[dotted,line width=0.75pt] coordinates { % h^3
    (1,\vara) (2,\vara*0.5*0.5*0.5) (3,\vara*0.25*0.25*0.25) (4,\vara*0.125*0.125*0.125) (5,\vara*0.0625*0.0625*0.0625)(6,\vara*0.03125*0.03125*0.03125)
    };
    \legend{$e_{L^2}^{\bu}$, $e_{PL^2}^{\bu}$, $e_{H^1}^{\bu}$, $e_{A}^{\bu}$, $e_{M}^{p}$, $\mathcal{O}(h^{2})$, $\mathcal{O}(h^{3})$}
  \end{semilogyaxis}
  \end{tikzpicture}
  }
%    \caption{Error norms for the sphere and $k=l=1$ with $\rho=\tilde\rho=1$.}
% \end{figure}    
\end{minipage}
\begin{minipage}{0.51\textwidth}
\scriptsize
      % \begin{figure}  
\scalebox{0.81}{
  \begin{tikzpicture}
  \def\vara{0.1}
  \def\varb{1.9}
  \def\varc{0.1}
\begin{semilogyaxis}[ xlabel={Refinement level}, ylabel={Error}, ymin=5E-9, ymax=5, legend style={ at={(0.5,1.02)}, anchor=south, legend columns =2, transpose legend }, legend cell align=left, cycle list name=mark list, title= {$\vect P_3$-$P_{2}$}, title style={at={(0.5,1.16)}, align=center, font=\Large} ]
    \addplot table[x=level, y=L2] {Ph1P3P2opt.dat};
    \addplot table[x=level, y=L2P] {Ph1P3P2opt.dat};
    \addplot table[x=level, y=realH1semi] {Ph1P3P2opt.dat};
    \addplot table[x=level, y=uerror] {Ph1P3P2opt.dat};
    \addplot table[x=level, y=merr] {Ph1P3P2opt.dat};
%    \addplot[dashed,line width=0.75pt] coordinates { % h
%    (1,\vara) (2,\vara*0.5) (3,\vara*0.25) (4,\vara*0.125) (5,\vara*0.0625)(6,\vara*0.03125)
%    };
%    \addplot[dashed,line width=0.75pt] coordinates { % h^2
%      (1,\varb) (2,\varb*0.5*0.5) (3,\varb*0.25*0.25) (4,\varb*0.125*0.125) (5,\varb*0.0625*0.0625)(6,\varb*0.03125*0.03125)
%    };
    \addplot[dotted,line width=0.75pt] coordinates { % h^3
    (1,\varb) (2,\varb*0.5*0.5*0.5) (3,\varb*0.25*0.25*0.25) (4,\varb*0.125*0.125*0.125) (5,\varb*0.0625*0.0625*0.0625)(6,\varb*0.03125*0.03125*0.03125)
    };
        \addplot[dashed,line width=0.75pt] coordinates { % h^4
    (1,\vara) (2,\vara*0.5*0.5*0.5*0.5) (3,\vara*0.25*0.25*0.25*0.25) (4,\vara*0.125*0.125*0.125*0.125) (5,\vara*0.0625*0.0625*0.0625*0.0625)(6,\vara*0.03125*0.03125*0.03125*0.03125)
    };
    \legend{$e_{L^2}^{\bu}$, $e_{PL^2}^{\bu}$, $e_{H^1}^{\bu}$, $e_{A}^{\bu}$, $e_{M}^{p}$, $\mathcal{O}(h^{3})$, $\mathcal{O}(h^{4})$}
  \end{semilogyaxis}
  \end{tikzpicture}
  }
%    \caption{Error norms for the sphere and $k=l=1$ with $\rho=\tilde\rho=1$.}
% \end{figure}       
\end{minipage}
\caption{Consistent formulation \eqref{discreteform1} on the unit sphere}\label{tablePh1sphereopt}
\end{figure}
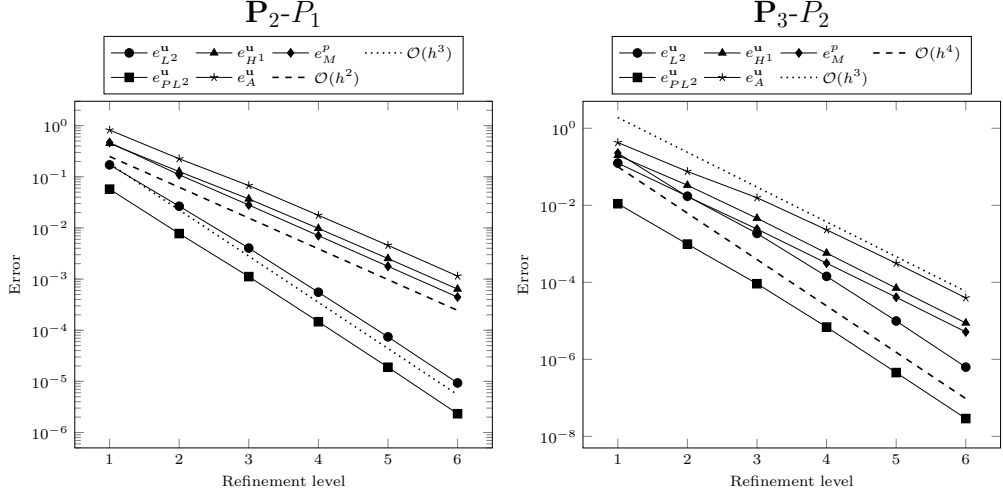

We begin with the consistent formulation \eqref{discreteform1}. In Figure \ref{tablePh1sphereopt} we show the errors for $\vect P_2$-$P_1$ and $\vect P_3$-$P_2$ finite elements.  We clearly observe optimal orders of convergence:  $e_{A}^{\bu} \sim h^k$, $e_{H^1}^{\bu} \sim  h^k$, $e_{M}^{p}\sim h^k$ and   $e_{L^2}^{\bu} \sim h^{k+1}$, $e_{PL^2}^{\bu}\sim h^{k+1}$. 

Concerning  the choice of the stabilization parameters $\rho_u$ and $\rho_p$ we note the following  (results of the experiments are not shown). For $\vect P_2$-$P_1$ finite elements and $\rho_u = h$, instead of $\rho_u = h^{-1}$, (and other parameters the same as above) we observe slightly slower than  $\mathcal{O}(h^2)$-convergence for the energy norm error $e_{A}^{\bu}$ and less than $\mathcal{O}(h^3)$-convergence for the $L^2$-error $e_{L^2}^{\bu}$. 
%The tangential $L^2$-error $e_{PL^2}^{\bu}$ is still of clear $\mathcal{O}(h^3)$-convergence. 
These suboptimal convergence orders are probably due to the consistency error (geometry error),  since taking superparametric finite elements, i.e. $k_g = 3$, leads to  optimal convergence orders.
%Taking into account the analysis of the vector-Laplace problem in \cite{jankuhn2019}, where we were only able to show optimal error bounds for $\rho_u = h^{-1}$, these results suggest that for the surface Stokes equation optimal error bounds only exist for $\rho_u = h^{-1}$. 
If we take $\rho_p = h^{-1}$ instead of $h$ (and other parameters the same as above), we observe a loss of one order in the  errors $e_{A}^{\bu}$ and $e_{M}^{p}$ and even a loss of one and a half order in  $e_{L^2}^{\bu}$. Taking $\rho_p = 1$ results in a loss of  a half order for  $e_{M}^{p}$, a loss of a quarter order for  $e_{L^2}^{\bu}$ and a loss of a half order for  $e_{PL^2}^{\bu}$.

\begin{figure}
\hspace*{-0.7cm}
\begin{minipage}{0.51\textwidth}
\scriptsize
      % \begin{figure}  
\scalebox{0.81}{
  \begin{tikzpicture}
  \def\vara{2.0}
  \def\varb{0.25}
  \def\varc{3.6}
\begin{semilogyaxis}[ xlabel={Refinement level}, ylabel={Error}, ymin=9E-7, ymax=8, legend style={ at={(0.5,1.02)}, anchor=south, legend columns =2, transpose legend }, legend cell align=left, cycle list name=mark list, title= {$\vect P_2$-$P_{1}$}, title style={at={(0.5,1.16)}, align=center, font=\Large} ]
    \addplot table[x=level, y=L2] {Ph2P2P1opt.dat};
    \addplot table[x=level, y=L2P] {Ph2P2P1opt.dat};
    \addplot table[x=level, y=realH1semi] {Ph2P2P1opt.dat};
    \addplot table[x=level, y=uerror] {Ph2P2P1opt.dat};
    \addplot table[x=level, y=merr] {Ph2P2P1opt.dat};
%    \addplot[dashed,line width=0.75pt] coordinates { % h
%    (1,\vara) (2,\vara*0.5) (3,\vara*0.25) (4,\vara*0.125) (5,\vara*0.0625)(6,\vara*0.03125)
%    };
    \addplot[densely dotted,line width=0.75pt] coordinates { % h^1.5
      (1,\varc) (2,\varc*0.5*0.707106781) (3,\varc*0.25*0.5) (4,\varc*0.125*0.353553391) (5,\varc*0.0625*0.25)(6,\varc*0.03125*0.176776695)
    };
    \addplot[dashed,line width=0.75pt] coordinates { % h^2
      (1,\vara) (2,\vara*0.5*0.5) (3,\vara*0.25*0.25) (4,\vara*0.125*0.125) (5,\vara*0.0625*0.0625)(6,\vara*0.03125*0.03125)
    };
    \addplot[loosely dotted,line width=0.75pt] coordinates { % h^3
    (1,\varb) (2,\varb*0.5*0.5*0.5) (3,\varb*0.25*0.25*0.25) (4,\varb*0.125*0.125*0.125) (5,\varb*0.0625*0.0625*0.0625)(6,\varb*0.03125*0.03125*0.03125)
    };
    \legend{$e_{L^2}^{\bu}$, $e_{PL^2}^{\bu}$, $e_{H^1}^{\bu}$, $e_{A}^{\bu}$, $e_{M}^{p}$, $\mathcal{O}(h^{1.5})$, $\mathcal{O}(h^{2})$, $\mathcal{O}(h^{3})$}
  \end{semilogyaxis}
  \end{tikzpicture}
  }
%    \caption{Error norms for the sphere and $k=l=1$ with $\rho=\tilde\rho=1$.}
% \end{figure}    
\end{minipage}
\begin{minipage}{0.51\textwidth}
\scriptsize
      % \begin{figure}  
\scalebox{0.81}{
  \begin{tikzpicture}
  \def\vara{0.8}
  \def\varb{0.5}
  \def\varc{0.15}
\begin{semilogyaxis}[ xlabel={Refinement level}, ylabel={Error}, ymin=7E-9, ymax=6, legend style={ at={(0.5,1.02)}, anchor=south, legend columns =2, transpose legend }, legend cell align=left, cycle list name=mark list, title= {$\vect P_3$-$P_{2}$}, title style={at={(0.5,1.16)}, align=center, font=\Large} ]
    \addplot table[x=level, y=L2] {Ph2P3P2opt.dat};
    \addplot table[x=level, y=L2P] {Ph2P3P2opt.dat};
    \addplot table[x=level, y=realH1semi] {Ph2P3P2opt.dat};
    \addplot table[x=level, y=uerror] {Ph2P3P2opt.dat};
    \addplot table[x=level, y=merr] {Ph2P3P2opt.dat};
%    \addplot[dashed,line width=0.75pt] coordinates { % h
%    (1,\varc) (2,\varc*0.5) (3,\varc*0.25) (4,\varc*0.125) (5,\varc*0.0625)(6,\varc*0.03125)
%    };
    \addplot[densely dotted,line width=0.75pt] coordinates { % h^2
      (1,\vara) (2,\vara*0.5*0.5) (3,\vara*0.25*0.25) (4,\vara*0.125*0.125) (5,\vara*0.0625*0.0625)(6,\vara*0.03125*0.03125)
    };
    \addplot[dashed,line width=0.75pt] coordinates { % h^3
    (1,\varb) (2,\varb*0.5*0.5*0.5) (3,\varb*0.25*0.25*0.25) (4,\varb*0.125*0.125*0.125) (5,\varb*0.0625*0.0625*0.0625)(6,\varb*0.03125*0.03125*0.03125)
    };
    \addplot[loosely dotted,line width=0.75pt] coordinates { % h^4
    (1,\varc) (2,\varc*0.5*0.5*0.5*0.5) (3,\varc*0.25*0.25*0.25*0.25) (4,\varc*0.125*0.125*0.125*0.125) (5,\varc*0.0625*0.0625*0.0625*0.0625)(6,\varc*0.03125*0.03125*0.03125*0.03125)
    };
    \legend{$e_{L^2}^{\bu}$, $e_{PL^2}^{\bu}$, $e_{H^1}^{\bu}$, $e_{A}^{\bu}$, $e_{M}^{p}$, $\mathcal{O}(h^{2})$, $\mathcal{O}(h^{3})$, $\mathcal{O}(h^{4})$}
  \end{semilogyaxis}
  \end{tikzpicture}
  }
%    \caption{Error norms for the sphere and $k=l=1$ with $\rho=\tilde\rho=1$.}
% \end{figure}       
\end{minipage}
\caption{Inconsistent formulation \eqref{discreteform2} on the unit sphere}\label{tablePh2sphereopt}
\end{figure}
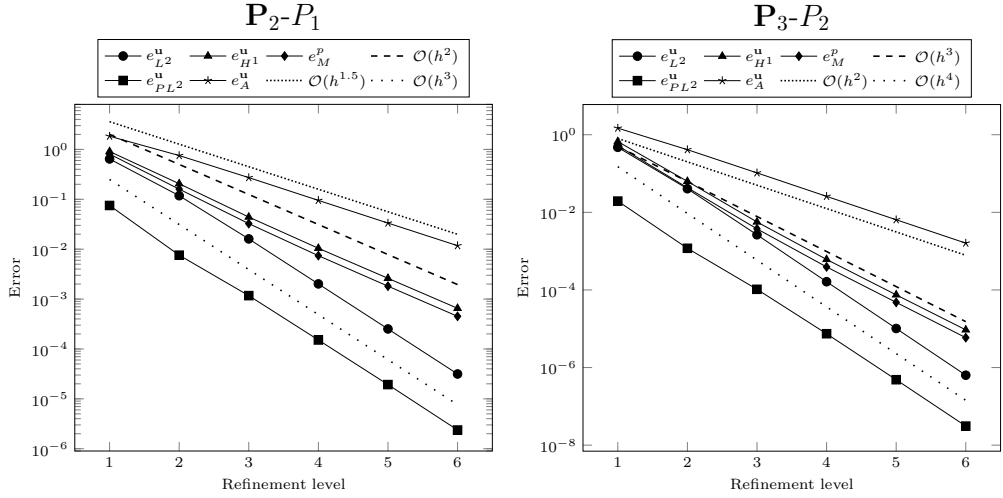

We now consider the inconsistent formulation \eqref{discreteform2}. In Figure \ref{tablePh2sphereopt} we show the errors for $\vect P_2$-$P_1$ and $\vect P_3$-$P_2$ finite elements. We observe $\mathcal{O}(h^{\frac12(k+1)})$-convergence for the energy norm error $e_{A}^{\bu}$, which is what we expect to see based on the analysis in \cite{jankuhn2019}. For the $e_{H^1}^{\bu}$- and $e_{M}^{p}$-errors we have $\mathcal{O}(h^{k})$-convergence and for the $L^2$-errors $e_{L^2}^{\bu}$ and $e_{PL^2}^{\bu}$ we see $\mathcal{O}(h^{k+1})$-convergence, which are all optimal.

A different scaling of the stabilization parameter $\rho_u$ does not have the same effect on the convergence behavior as described above for the consistent method \eqref{discreteform1}. For $\vect P_2$-$P_1$ finite elements and $\rho_u = h$ we still observe the same optimal convergence order as for $\rho_u = h^{-1}$. 
For the stabilization parameter $\rho_p$, however,  we see similar effects as described above for the consistent formulation \eqref{discreteform1}. \\

\begin{figure}
\hspace*{-0.7cm}
\begin{minipage}{0.51\textwidth}
\scriptsize
      % \begin{figure}  
\scalebox{0.81}{
  \begin{tikzpicture}
  \def\vara{0.18}
  \def\varb{0.25}
  \def\varc{0.1}
\begin{semilogyaxis}[ xlabel={Refinement level}, ylabel={Error}, ymin=2E-6, ymax=3, legend style={ at={(0.5,1.02)}, anchor=south, legend columns =2, transpose legend }, legend cell align=left, cycle list name=mark list, title= {$\vect P_2$-$P_{1}$}, title style={at={(0.5,1.16)}, align=center, font=\Large} ]
    \addplot table[x=level, y=L2P1h] {P2P1comp.dat};
    \addplot table[x=level, y=realH1semiP1h] {P2P1comp.dat};
        \addplot table[x=level, y=L2P2h] {P2P1comp.dat};
    \addplot table[x=level, y=realH1semiP2h] {P2P1comp.dat};
%    \addplot[dashed,line width=0.75pt] coordinates { % h
%    (1,\vara) (2,\vara*0.5) (3,\vara*0.25) (4,\vara*0.125) (5,\vara*0.0625)(6,\vara*0.03125)
%    };
    \addplot[dashed,line width=0.75pt] coordinates { % h^2
      (1,\varb) (2,\varb*0.5*0.5) (3,\varb*0.25*0.25) (4,\varb*0.125*0.125) (5,\varb*0.0625*0.0625)(6,\varb*0.03125*0.03125)
    };
    \addplot[dotted,line width=0.75pt] coordinates { % h^3
    (1,\vara) (2,\vara*0.5*0.5*0.5) (3,\vara*0.25*0.25*0.25) (4,\vara*0.125*0.125*0.125) (5,\vara*0.0625*0.0625*0.0625)(6,\vara*0.03125*0.03125*0.03125)
    };
    \legend{$e_{L^2}^{\bu}$ \eqref{discreteform1}, $e_{H^1}^{\bu}$ \eqref{discreteform1}, $e_{L^2}^{\bu}$ \eqref{discreteform2}, $e_{H^1}^{\bu}$ \eqref{discreteform2}, $\mathcal{O}(h^{2})$, $\mathcal{O}(h^{3})$}
  \end{semilogyaxis}
  \end{tikzpicture}
  }
%    \caption{Error norms for the sphere and $k=l=1$ with $\rho=\tilde\rho=1$.}
% \end{figure}    
\end{minipage}
\begin{minipage}{0.51\textwidth}
\scriptsize
      % \begin{figure}  
\scalebox{0.81}{
  \begin{tikzpicture}
  \def\vara{0.1}
  \def\varb{1.2}
  \def\varc{0.1}
\begin{semilogyaxis}[ xlabel={Refinement level}, ylabel={Error}, ymin=2E-8, ymax=5, legend style={ at={(0.5,1.02)}, anchor=south, legend columns =2, transpose legend }, legend cell align=left, cycle list name=mark list, title= {$\vect P_3$-$P_{2}$}, title style={at={(0.5,1.16)}, align=center, font=\Large} ]
    \addplot table[x=level, y=L2P1h] {P3P2comp.dat};
    \addplot table[x=level, y=realH1semiP1h] {P3P2comp.dat};
        \addplot table[x=level, y=L2P2h] {P3P2comp.dat};
    \addplot table[x=level, y=realH1semiP2h] {P3P2comp.dat};
%    \addplot[dashed,line width=0.75pt] coordinates { % h
%    (1,\vara) (2,\vara*0.5) (3,\vara*0.25) (4,\vara*0.125) (5,\vara*0.0625)(6,\vara*0.03125)
%    };
%    \addplot[dashed,line width=0.75pt] coordinates { % h^2
%      (1,\varb) (2,\varb*0.5*0.5) (3,\varb*0.25*0.25) (4,\varb*0.125*0.125) (5,\varb*0.0625*0.0625)(6,\varb*0.03125*0.03125)
%    };
    \addplot[dotted,line width=0.75pt] coordinates { % h^3
    (1,\varb) (2,\varb*0.5*0.5*0.5) (3,\varb*0.25*0.25*0.25) (4,\varb*0.125*0.125*0.125) (5,\varb*0.0625*0.0625*0.0625)(6,\varb*0.03125*0.03125*0.03125)
    };
        \addplot[dashed,line width=0.75pt] coordinates { % h^4
    (1,\vara) (2,\vara*0.5*0.5*0.5*0.5) (3,\vara*0.25*0.25*0.25*0.25) (4,\vara*0.125*0.125*0.125*0.125) (5,\vara*0.0625*0.0625*0.0625*0.0625)(6,\vara*0.03125*0.03125*0.03125*0.03125)
    };
    \legend{$e_{L^2}^{\bu}$ \eqref{discreteform1}, $e_{H^1}^{\bu}$ \eqref{discreteform1}, $e_{L^2}^{\bu}$ \eqref{discreteform2}, $e_{H^1}^{\bu}$ \eqref{discreteform2}, $\mathcal{O}(h^{3})$, $\mathcal{O}(h^{4})$}
  \end{semilogyaxis}
  \end{tikzpicture}
  }
%    \caption{Error norms for the sphere and $k=l=1$ with $\rho=\tilde\rho=1$.}
% \end{figure}       
\end{minipage}
\caption{Comparison of \eqref{discreteform1} and \eqref{discreteform2} on the unit sphere}\label{tablecompsphere}
\end{figure}
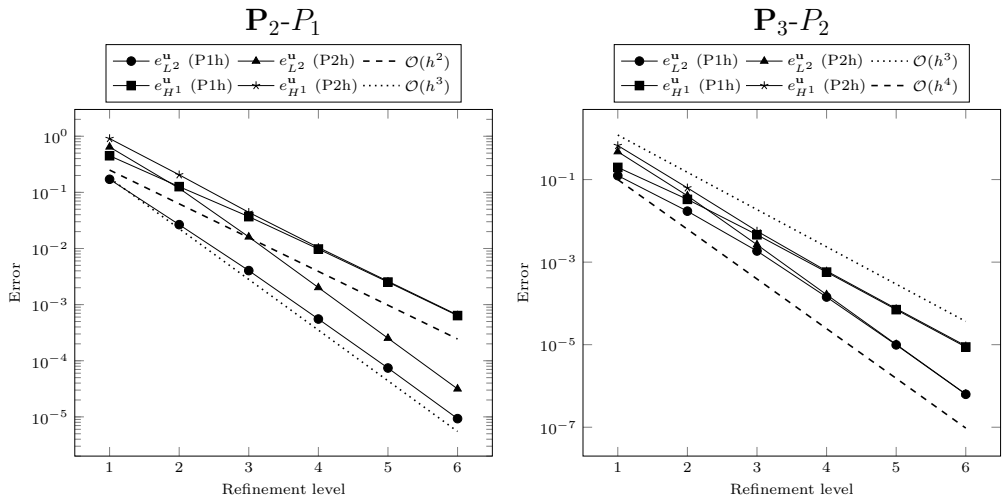

 Both methods \eqref{discreteform1} and \eqref{discreteform2} have optimal order errors $e_{L^2}^{\bu}$ and $e_{H^1}^{\bu}$. The question arises which of the two methods results in a smaller absolute error. Therefore, in Figure \ref{tablecompsphere} we show the $e_{L^2}^{\bu}$- and $e_{H^1}^{\bu}$-errors for both methods in one plot. For $\vect P_2$-$P_1$ finite elements the $e_{H^1}^{\bu}$-errors  differ only slightly for the first three refinement levels and the $e_{L^2}^{\bu}$-error of the consistent formulation is approximately one order of magnitude smaller than the one of the inconsistent formulation. For $\vect P_3$-$P_2$ finite elements the errors of both methods are almost the same. 
% For both cases of finite elements the $e_{L^2}^{\bu}$- and $e_{H^1}^{\bu}$-errors are smaller for the consistent formulation \eqref{discreteform1}.

\begin{remark} \rm
A special situation occurs if one considers a Stokes problem on the sphere with a divergence-free velocity solution $\bu$. In such a case the energy norm of the inconsistent method is not of order $e_{A}^{\bu} \sim h^{\frac12(k+1)}$ (as in the results above), but of order $e_{A}^{\bu} \sim h^{k}$. This improvement can be explained as follows. From the analysis in \cite{jankuhn2019} we notice that for the inconsistent method the dominant  inconsistency term is
\begin{equation*}
 ( E(\bu) , (\bv_h^l \cdot \bn)\bH)_{L^2(\Gamma)} = \int_{\Gamma} (\bv_h^l \cdot \bn) \textrm{tr} (E(\bu)  \bH) \, ds,
\end{equation*}
with $\bv_h^l$ the lifting of a finite element function from $\Gamma_h$ to $\Gamma$. 
For the sphere we have $\bH = \bP$ and thus 
\begin{equation*}
\textrm{tr} (E(\bu)  \bH) = \divG(\bu),   
\end{equation*}
which vanishes for a divergence-free solution $\bu$.
\end{remark}

\subsubsection{Effects related to the penalty term} \label{sectionpenaltyterm}

As mentioned in Section~\ref{sectionanalysis}, in case of the vector-Laplace problem, for optimal convergence it is essential that one uses $k_p=k+1$ (i.e. a one order better approximation for the normal approximation $\tilde \bn_h$). For the Stokes problem we performed an experiment with $\vect P_2$-$P_1$ finite elements in which all parameters and components are the same as in the experiments above, except for $k_p$: we take $k_p=2$ instead of $k_p=3$. The results are presented in Figure \ref{tablekpsphere}. In case of the consistent formulation \eqref{discreteform1} we lose, compared to $k_p=3$, one order for the $L^2$-errors $e_{L^2}^{\bu}$ and $e_{PL^2}^{\bu}$ and one order for the energy norm error $e_{A}^{\bu}$. The convergence of the $e_{H^1}^{\bu}$-error and $e_{M}^{p}$-error is a little worse than $\mathcal{O}(h^{2})$. For the inconsistent formulation \eqref{discreteform2} the effect is even stronger. In that case the $L^2$-errors $e_{L^2}^{\bu}$ and $e_{PL^2}^{\bu}$ are only of order $\mathcal{O}(
h)$ and the energy norm error $e_{A}^{\bu}$ converges significantly slower than first order. The convergence of the $e_{H^1}^{\bu}$-error and $e_{M}^{p}$-error is a little worse than $\mathcal{O}(h^{2})$.
%The $e_{H^1}^{\bu}$-error and the energy norm error $e_{M}^{p}$ are also of a little less than $\mathcal{O}(h^{2})$-convergence order.

\begin{figure}
\hspace*{-0.7cm}
\begin{minipage}{0.51\textwidth}
\scriptsize
      % \begin{figure}  
\scalebox{0.81}{
  \begin{tikzpicture}
  \def\vara{0.25}
  \def\varb{0.25}
  \def\varc{0.1}
\begin{semilogyaxis}[ xlabel={Refinement level}, ylabel={Error}, ymin=2E-6, ymax=3, legend style={ at={(0.5,1.02)}, anchor=south, legend columns =2, transpose legend }, legend cell align=left, cycle list name=mark list, title= {\eqref{discreteform1}}, title style={at={(0.5,1.16)}, align=center, font=\Large} ]
    \addplot table[x=level, y=L2] {Ph1kp.dat};
    \addplot table[x=level, y=L2P] {Ph1kp.dat};
    \addplot table[x=level, y=realH1semi] {Ph1kp.dat};
    \addplot table[x=level, y=uerror] {Ph1kp.dat};
    \addplot table[x=level, y=merr] {Ph1kp.dat};
    \addplot[dotted,line width=0.75pt] coordinates { % h
    (1,\vara) (2,\vara*0.5) (3,\vara*0.25) (4,\vara*0.125) (5,\vara*0.0625)(6,\vara*0.03125)
    };
    \addplot[dashed,line width=0.75pt] coordinates { % h^2
      (1,\varb) (2,\varb*0.5*0.5) (3,\varb*0.25*0.25) (4,\varb*0.125*0.125) (5,\varb*0.0625*0.0625)(6,\varb*0.03125*0.03125)
    };
%    \addplot[dotted,line width=0.75pt] coordinates { % h^3
%    (1,\vara) (2,\vara*0.5*0.5*0.5) (3,\vara*0.25*0.25*0.25) (4,\vara*0.125*0.125*0.125) (5,\vara*0.0625*0.0625*0.0625)(6,\vara*0.03125*0.03125*0.03125)
%    };
    \legend{$e_{L^2}^{\bu}$, $e_{PL^2}^{\bu}$, $e_{H^1}^{\bu}$, $e_{A}^{\bu}$, $e_{M}^{p}$, $\mathcal{O}(h)$, $\mathcal{O}(h^{2})$}
  \end{semilogyaxis}
  \end{tikzpicture}
  }
%    \caption{Error norms for the sphere and $k=l=1$ with $\rho=\tilde\rho=1$.}
% \end{figure}    
\end{minipage}
\begin{minipage}{0.51\textwidth}
\scriptsize
      % \begin{figure}  
\scalebox{0.81}{
  \begin{tikzpicture}
  \def\vara{1.2}
  \def\varb{0.25}
  \def\varc{0.005}
\begin{semilogyaxis}[ xlabel={Refinement level}, ylabel={Error}, ymin=8E-5, ymax=4, legend style={ at={(0.5,1.02)}, anchor=south, legend columns =2, transpose legend }, legend cell align=left, cycle list name=mark list, title= {\eqref{discreteform2}}, title style={at={(0.5,1.16)}, align=center, font=\Large} ]
    \addplot table[x=level, y=L2] {Ph2kp.dat};
    \addplot table[x=level, y=L2P] {Ph2kp.dat};
    \addplot table[x=level, y=realH1semi] {Ph2kp.dat};
    \addplot table[x=level, y=uerror] {Ph2kp.dat};
    \addplot table[x=level, y=merr] {Ph2kp.dat};
    \addplot[dotted,line width=0.75pt] coordinates { % h
    (1,\varc) (2,\varc*0.5) (3,\varc*0.25) (4,\varc*0.125) (5,\varc*0.0625)(6,\varc*0.03125)
    };
    \addplot[dashed,line width=0.75pt] coordinates { % h^2
      (1,\vara) (2,\vara*0.5*0.5) (3,\vara*0.25*0.25) (4,\vara*0.125*0.125) (5,\vara*0.0625*0.0625)(6,\vara*0.03125*0.03125)
    };
%    \addplot[dotted,line width=0.75pt] coordinates { % h^3
%    (1,\varb) (2,\varb*0.5*0.5*0.5) (3,\varb*0.25*0.25*0.25) (4,\varb*0.125*0.125*0.125) (5,\varb*0.0625*0.0625*0.0625)(6,\varb*0.03125*0.03125*0.03125)
%    };
    \legend{$e_{L^2}^{\bu}$, $e_{PL^2}^{\bu}$, $e_{H^1}^{\bu}$, $e_{A}^{\bu}$, $e_{M}^{p}$, $\mathcal{O}(h)$, $\mathcal{O}(h^{2})$}
  \end{semilogyaxis}
  \end{tikzpicture}
  }
%    \caption{Error norms for the sphere and $k=l=1$ with $\rho=\tilde\rho=1$.}
% \end{figure}       
\end{minipage}
\caption{$\vect P_2$-$P_1$ finite elements with $k_p=2$}\label{tablekpsphere}
\end{figure} 

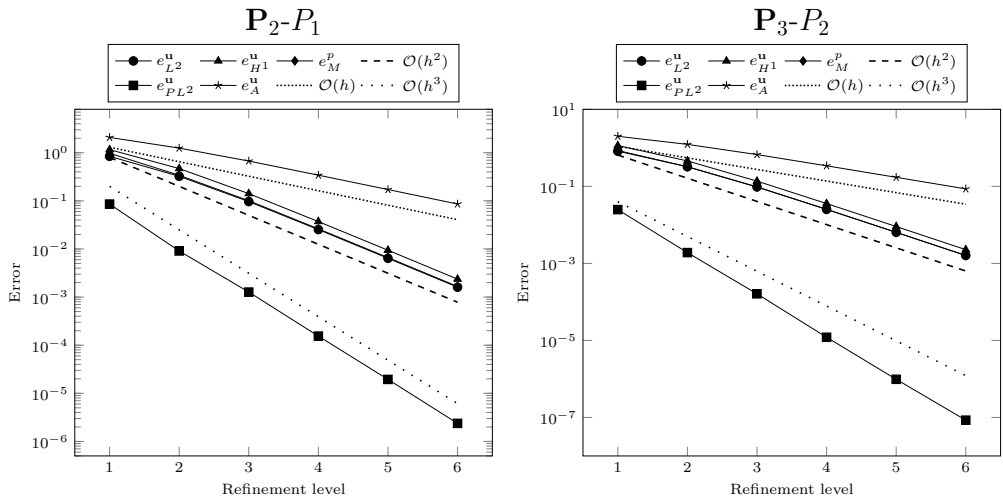
\begin{figure}
\hspace*{-0.7cm}
\begin{minipage}{0.51\textwidth}
\scriptsize
      % \begin{figure}  
\scalebox{0.81}{
  \begin{tikzpicture}
  \def\vara{0.8}
  \def\varb{0.2}
  \def\varc{1.3}
\begin{semilogyaxis}[ xlabel={Refinement level}, ylabel={Error}, ymin=5E-7, ymax=8, legend style={ at={(0.5,1.02)}, anchor=south, legend columns =2, transpose legend }, legend cell align=left, cycle list name=mark list, title= {$\vect P_2$-$P_{1}$}, title style={at={(0.5,1.16)}, align=center, font=\Large} ]
    \addplot table[x=level, y=L2] {Ph2P2P1eta.dat};
    \addplot table[x=level, y=L2P] {Ph2P2P1eta.dat};
    \addplot table[x=level, y=realH1semi] {Ph2P2P1eta.dat};
    \addplot table[x=level, y=uerror] {Ph2P2P1eta.dat};
    \addplot table[x=level, y=merr] {Ph2P2P1eta.dat};
    \addplot[densely dotted,line width=0.75pt] coordinates { % h
    (1,\varc) (2,\varc*0.5) (3,\varc*0.25) (4,\varc*0.125) (5,\varc*0.0625)(6,\varc*0.03125)
    };
    \addplot[dashed,line width=0.75pt] coordinates { % h^2
      (1,\vara) (2,\vara*0.5*0.5) (3,\vara*0.25*0.25) (4,\vara*0.125*0.125) (5,\vara*0.0625*0.0625)(6,\vara*0.03125*0.03125)
    };
    \addplot[ loosely dotted,line width=0.75pt] coordinates { % h^3
    (1,\varb) (2,\varb*0.5*0.5*0.5) (3,\varb*0.25*0.25*0.25) (4,\varb*0.125*0.125*0.125) (5,\varb*0.0625*0.0625*0.0625)(6,\varb*0.03125*0.03125*0.03125)
    };
    \legend{$e_{L^2}^{\bu}$, $e_{PL^2}^{\bu}$, $e_{H^1}^{\bu}$, $e_{A}^{\bu}$, $e_{M}^{p}$, $\mathcal{O}(h)$, $\mathcal{O}(h^{2})$, $\mathcal{O}(h^{3})$}
  \end{semilogyaxis}
  \end{tikzpicture}
  }
%    \caption{Error norms for the sphere and $k=l=1$ with $\rho=\tilde\rho=1$.}
% \end{figure}    
\end{minipage}
\begin{minipage}{0.51\textwidth}
\scriptsize
      % \begin{figure}  
\scalebox{0.81}{
  \begin{tikzpicture}
  \def\vara{0.65}
  \def\varb{0.04}
  \def\varc{1.1}
\begin{semilogyaxis}[ xlabel={Refinement level}, ylabel={Error}, ymin=1E-8, ymax=10, legend style={ at={(0.5,1.02)}, anchor=south, legend columns =2, transpose legend }, legend cell align=left, cycle list name=mark list, title= {$\vect P_3$-$P_{2}$}, title style={at={(0.5,1.16)}, align=center, font=\Large} ]
    \addplot table[x=level, y=L2] {Ph2P3P2eta.dat};
    \addplot table[x=level, y=L2P] {Ph2P3P2eta.dat};
    \addplot table[x=level, y=realH1semi] {Ph2P3P2eta.dat};
    \addplot table[x=level, y=uerror] {Ph2P3P2eta.dat};
    \addplot table[x=level, y=merr] {Ph2P3P2eta.dat};
    \addplot[densely dotted,line width=0.75pt] coordinates { % h
    (1,\varc) (2,\varc*0.5) (3,\varc*0.25) (4,\varc*0.125) (5,\varc*0.0625)(6,\varc*0.03125)
    };
    \addplot[dashed,line width=0.75pt] coordinates { % h^2
      (1,\vara) (2,\vara*0.5*0.5) (3,\vara*0.25*0.25) (4,\vara*0.125*0.125) (5,\vara*0.0625*0.0625)(6,\vara*0.03125*0.03125)
    };
    \addplot[loosely dotted,line width=0.75pt] coordinates { % h^3
    (1,\varb) (2,\varb*0.5*0.5*0.5) (3,\varb*0.25*0.25*0.25) (4,\varb*0.125*0.125*0.125) (5,\varb*0.0625*0.0625*0.0625)(6,\varb*0.03125*0.03125*0.03125)
    };
    \legend{$e_{L^2}^{\bu}$, $e_{PL^2}^{\bu}$, $e_{H^1}^{\bu}$, $e_{A}^{\bu}$, $e_{M}^{p}$, $\mathcal{O}(h)$, $\mathcal{O}(h^{2})$, $\mathcal{O}(h^{3})$}
  \end{semilogyaxis}
  \end{tikzpicture}
  }
%    \caption{Error norms for the sphere and $k=l=1$ with $\rho=\tilde\rho=1$.}
% \end{figure}       
\end{minipage}
\caption{Inconsistent formulation \eqref{discreteform2} with $\eta=h^{-2}$}\label{tableetasphere}
\end{figure}

As noted in Section \ref{sectionanalysis}, to obtain optimal convergence for the inconsistent formulation of the vector-Laplace problem the scaling of the penalty parameter $\eta$ has to depend on the degree of the finite element space $k$: $\eta \sim  h^{-(k+1)}$.  For the inconsistent formulation of the Stokes problem \eqref{discreteform2} we performed an experiment with $\vect P_2$-$P_1$ and $\vect P_3$-$P_2$ finite elements in which all parameters and components are the same as in the experiments above, except for $\eta$: we take $\eta=h^{-2}$ instead of $\eta = h^{-(k+1)}$. The results are shown in Figure \ref{tableetasphere}. For $\vect P_2$-$P_1$ finite elements we observe suboptimal $\mathcal{O}(h)$-convergence for the energy norm error $e_{A}^{\bu}$, which is half an order less than for $\eta = h^{-3}$. We still have $\mathcal{O}(h^{2})$-convergence for the $e_{H^1}^{\bu}$- and $e_{M}^{p}$-errors and $\mathcal{O}(h^{3})$-convergence for the tangential $L^2$-error $e_{PL^2}^{\bu}$, which are both 
optimal. The full $L^2$-error $e_{L^2}^{\bu}$ though loses one order compared to $\eta = h^{-3}$. For $\vect P_3$-$P_2$ finite elements we observe the same convergence orders for all the errors except for the tangential $L^2$-error $e_{PL^2}^{\bu}$, which is a bit better as for the $\vect P_2$-$P_1$ case. Hence, for $\vect P_3$-$P_2$ finite elements all errors show suboptimal convergence if we take $\eta=h^{-2}$.

\begin{figure}
\hspace*{-0.7cm}
\begin{minipage}{0.51\textwidth}
\scriptsize
      % \begin{figure}  
\scalebox{0.81}{
  \begin{tikzpicture}
  \def\vara{2.0}
  \def\varb{0.25}
  \def\varc{3.6}
\begin{semilogyaxis}[ xlabel={Refinement level}, ylabel={Error}, ymin=2E-4, ymax=8, legend style={ at={(0.5,1.02)}, anchor=south, legend columns =2, transpose legend }, legend cell align=left, cycle list name=mark list, title= {$\vect P_2$-$P_{1}$}, title style={at={(0.5,1.19)}, align=center, font=\Large} ]
    \addplot table[x=level, y=ah] {Ph2P2P1exp.dat};
    \addplot table[x=level, y=sh] {Ph2P2P1exp.dat};
    \addplot table[x=level, y=kh] {Ph2P2P1exp.dat};
    \addplot table[x=level, y=uerror] {Ph2P2P1exp.dat};
%    \addplot[dashed,line width=0.75pt] coordinates { % h
%    (1,\vara) (2,\vara*0.5) (3,\vara*0.25) (4,\vara*0.125) (5,\vara*0.0625)(6,\vara*0.03125)
%    };
    \addplot[densely dotted,line width=0.75pt] coordinates { % h^1.5
      (1,\varc) (2,\varc*0.5*0.707106781) (3,\varc*0.25*0.5) (4,\varc*0.125*0.353553391) (5,\varc*0.0625*0.25)(6,\varc*0.03125*0.176776695)
    };
    \addplot[dashed,line width=0.75pt] coordinates { % h^2
      (1,\vara) (2,\vara*0.5*0.5) (3,\vara*0.25*0.25) (4,\vara*0.125*0.125) (5,\vara*0.0625*0.0625)(6,\vara*0.03125*0.03125)
    };
%    \addplot[loosely dotted,line width=0.75pt] coordinates { % h^3
%    (1,\varb) (2,\varb*0.5*0.5*0.5) (3,\varb*0.25*0.25*0.25) (4,\varb*0.125*0.125*0.125) (5,\varb*0.0625*0.0625*0.0625)(6,\varb*0.03125*0.03125*0.03125)
%    };
    \legend{$a_h({\mathbf{e}_h,\mathbf{e}_h})^{\frac12}$, $s_h({\mathbf{e}_h,\mathbf{e}_h})^{\frac12}$, $k_h({\mathbf{e}_h,\mathbf{e}_h})^{\frac12}$, $e_{A}^{\bu}$, $\mathcal{O}(h^{1.5})$, $\mathcal{O}(h^{2})$}
  \end{semilogyaxis}
  \end{tikzpicture}
  }
%    \caption{Error norms for the sphere and $k=l=1$ with $\rho=\tilde\rho=1$.}
% \end{figure}    
\end{minipage}
\begin{minipage}{0.51\textwidth}
\scriptsize
      % \begin{figure}  
\scalebox{0.81}{
  \begin{tikzpicture}
  \def\vara{0.8}
  \def\varb{0.5}
  \def\varc{0.15}
\begin{semilogyaxis}[ xlabel={Refinement level}, ylabel={Error}, ymin=2E-6, ymax=6, legend style={ at={(0.5,1.02)}, anchor=south, legend columns =2, transpose legend }, legend cell align=left, cycle list name=mark list, title= {$\vect P_3$-$P_{2}$}, title style={at={(0.5,1.19)}, align=center, font=\Large} ]
    \addplot table[x=level, y=ah] {Ph2P3P2exp.dat};
    \addplot table[x=level, y=sh] {Ph2P3P2exp.dat};
    \addplot table[x=level, y=kh] {Ph2P3P2exp.dat};
    \addplot table[x=level, y=uerror] {Ph2P3P2exp.dat};
%    \addplot[dashed,line width=0.75pt] coordinates { % h
%    (1,\varc) (2,\varc*0.5) (3,\varc*0.25) (4,\varc*0.125) (5,\varc*0.0625)(6,\varc*0.03125)
%    };
    \addplot[densely dotted,line width=0.75pt] coordinates { % h^2
      (1,\vara) (2,\vara*0.5*0.5) (3,\vara*0.25*0.25) (4,\vara*0.125*0.125) (5,\vara*0.0625*0.0625)(6,\vara*0.03125*0.03125)
    };
    \addplot[dashed,line width=0.75pt] coordinates { % h^3
    (1,\varb) (2,\varb*0.5*0.5*0.5) (3,\varb*0.25*0.25*0.25) (4,\varb*0.125*0.125*0.125) (5,\varb*0.0625*0.0625*0.0625)(6,\varb*0.03125*0.03125*0.03125)
    };
%    \addplot[loosely dotted,line width=0.75pt] coordinates { % h^4
%    (1,\varc) (2,\varc*0.5*0.5*0.5*0.5) (3,\varc*0.25*0.25*0.25*0.25) (4,\varc*0.125*0.125*0.125*0.125) (5,\varc*0.0625*0.0625*0.0625*0.0625)(6,\varc*0.03125*0.03125*0.03125*0.03125)
%    };
    \legend{$a_h({\mathbf{e}_h,\mathbf{e}_h})^{\frac12}$, $s_h({\mathbf{e}_h,\mathbf{e}_h})^{\frac12}$, $k_h({\mathbf{e}_h,\mathbf{e}_h})^{\frac12}$, $e_{A}^{\bu}$, $\mathcal{O}(h^{2})$, $\mathcal{O}(h^{3})$}
  \end{semilogyaxis}
  \end{tikzpicture}
  }
%    \caption{Error norms for the sphere and $k=l=1$ with $\rho=\tilde\rho=1$.}
% \end{figure}       
\end{minipage}
\caption{Components of the energy norm error for \eqref{discreteform2}}\label{tablePh2sphereexp}
\end{figure}
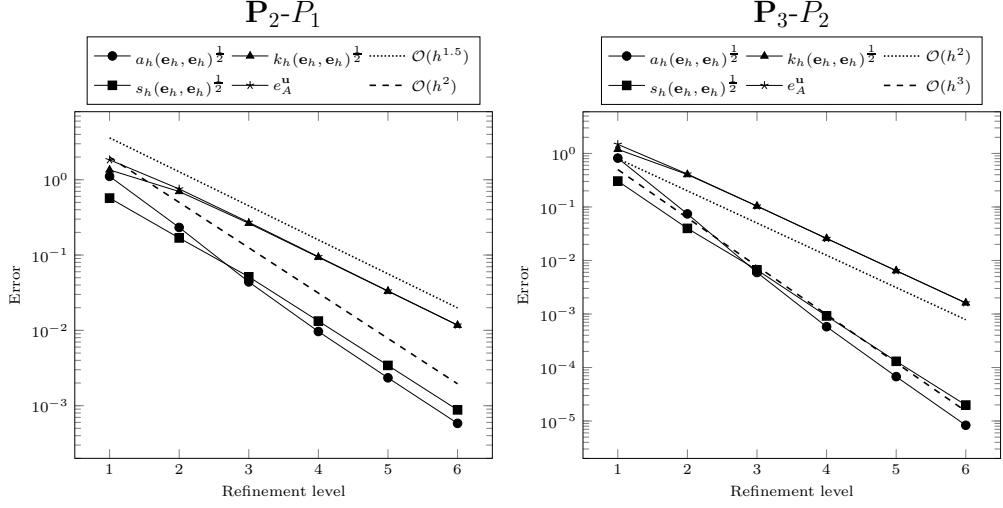

Finally we briefly discuss the $\mathcal{O}(h^{\frac12(k+1)})$-convergence in the energy norm error $e_A^{\bu}$ observed for the inconsistent method, cf. Figure~\ref{tablePh2sphereopt}. We call this convergence rate ``optimal'', due to the penalty term which is included in the energy norm:
\[
  \|\bv\|_{A}^2=A_{h}(\bv,\bv) := a_{h}(\bv,\bv) + s_h(\bv,\bv)+k_h(\bv,\bv).
\]
 For all parameters and components we take the default values. For the penalty term part of the energy norm error  $ \|\mathbf e_h\|_A=e_{A}^{\bu}$ we have
\[
  k_h(\mathbf e_h,\mathbf e_h)^\frac12= \eta^\frac12 \|\tilde \bn_h \cdot \mathbf e_h\|_{L^2(\Gamma_h)} = h^{-\frac12(k+1)}\|\tilde \bn_h \cdot \mathbf e_h\|_{L^2(\Gamma_h)}.
\]
For the term $\|\tilde \bn_h \cdot \mathbf e_h\|_{L^2(\Gamma_h)}$ the best one can expect (based on an interpolation error) is $\|\tilde \bn_h \cdot \mathbf e_h\|_{L^2(\Gamma_h)} \sim h^{k+1}$. Hence, for the penalty term part of the energy norm error we obtain an optimal convergence rate $ k_h(\mathbf e_h,\mathbf e_h)^\frac12 \sim h^{\frac12(k+1)}$.  This explains  why $h^{\frac12(k+1)}$ is the ``optimal'' convergence rate for the energy norm error $e_{A}^{\bu}$, which is indeed attained for the inconsistent method and $k=2$, $k=3$, cf. ~Figure~\ref{tablePh2sphereopt}.
To illustrate this, we  performed an experiment in which the three different contributions to the energy norm error are shown separately. 
In Figure \ref{tablePh2sphereexp}, for $\vect P_2$-$P_1$ and $\vect P_3$-$P_2$ finite elements, we show the energy norm $\|\mathbf e_h\|_A=e_{A}^{\bu}$   and its three components $a_{h}(\mathbf e_h,\mathbf e_h)^\frac12$,
$s_h(\mathbf e_h,\mathbf e_h)^\frac12$, $k_h(\mathbf e_h,\mathbf e_h)^\frac12$. We clearly observe that $e_A^{\bu} \approx k_h(\mathbf e_h,\mathbf e_h)^\frac12 \sim h^{\frac12(k+1)}$. Furthermore, the other error components have a higher rate of convergence: $a_h(\mathbf e_h,\mathbf e_h)^\frac12 \sim h^{k}$ and $s_h(\mathbf e_h,\mathbf e_h)^\frac12 \sim h^{k}$. 

\subsection{Results for the torus}
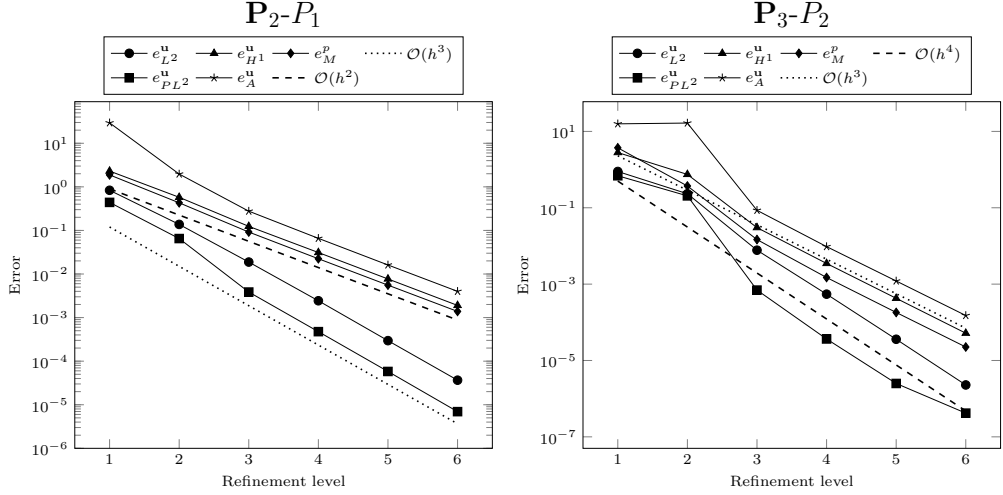
\begin{figure}
\hspace*{-0.7cm}
\begin{minipage}{0.51\textwidth}
\scriptsize
      % \begin{figure}  
\scalebox{0.81}{
  \begin{tikzpicture}
  \def\vara{0.9}
  \def\varb{0.12}
  \def\varc{0.1}
\begin{semilogyaxis}[ xlabel={Refinement level}, ylabel={Error}, ymin=1E-6, ymax=90, legend style={ at={(0.5,1.02)}, anchor=south, legend columns =2, transpose legend }, legend cell align=left, cycle list name=mark list, title= {$\vect P_2$-$P_{1}$}, title style={at={(0.5,1.16)}, align=center, font=\Large} ]
    \addplot table[x=level, y=L2] {Ph1P2P1torusopt.dat};
    \addplot table[x=level, y=L2P] {Ph1P2P1torusopt.dat};
    \addplot table[x=level, y=realH1semi] {Ph1P2P1torusopt.dat};
    \addplot table[x=level, y=uerror] {Ph1P2P1torusopt.dat};
    \addplot table[x=level, y=merr] {Ph1P2P1torusopt.dat};
%    \addplot[dashed,line width=0.75pt] coordinates { % h
%    (1,\vara) (2,\vara*0.5) (3,\vara*0.25) (4,\vara*0.125) (5,\vara*0.0625)(6,\vara*0.03125)
%    };
    \addplot[dashed,line width=0.75pt] coordinates { % h^2
      (1,\vara) (2,\vara*0.5*0.5) (3,\vara*0.25*0.25) (4,\vara*0.125*0.125) (5,\vara*0.0625*0.0625)(6,\vara*0.03125*0.03125)
    };
    \addplot[dotted,line width=0.75pt] coordinates { % h^3
    (1,\varb) (2,\varb*0.5*0.5*0.5) (3,\varb*0.25*0.25*0.25) (4,\varb*0.125*0.125*0.125) (5,\varb*0.0625*0.0625*0.0625)(6,\varb*0.03125*0.03125*0.03125)
    };
    \legend{$e_{L^2}^{\bu}$, $e_{PL^2}^{\bu}$, $e_{H^1}^{\bu}$, $e_{A}^{\bu}$, $e_{M}^{p}$, $\mathcal{O}(h^{2})$, $\mathcal{O}(h^{3})$}
  \end{semilogyaxis}
  \end{tikzpicture}
  }
%    \caption{Error norms for the sphere and $k=l=1$ with $\rho=\tilde\rho=1$.}
% \end{figure}    
\end{minipage}
\begin{minipage}{0.51\textwidth}
\scriptsize
      % \begin{figure}  
\scalebox{0.81}{
  \begin{tikzpicture}
  \def\vara{0.5}
  \def\varb{2.3}
  \def\varc{0.1}
\begin{semilogyaxis}[ xlabel={Refinement level}, ylabel={Error}, ymin=5E-8, ymax=60, legend style={ at={(0.5,1.02)}, anchor=south, legend columns =2, transpose legend }, legend cell align=left, cycle list name=mark list, title= {$\vect P_3$-$P_{2}$}, title style={at={(0.5,1.16)}, align=center, font=\Large} ]
    \addplot table[x=level, y=L2] {Ph1P3P2torusopt.dat};
    \addplot table[x=level, y=L2P] {Ph1P3P2torusopt.dat};
    \addplot table[x=level, y=realH1semi] {Ph1P3P2torusopt.dat};
    \addplot table[x=level, y=uerror] {Ph1P3P2torusopt.dat};
    \addplot table[x=level, y=merr] {Ph1P3P2torusopt.dat};
%    \addplot[dashed,line width=0.75pt] coordinates { % h
%    (1,\vara) (2,\vara*0.5) (3,\vara*0.25) (4,\vara*0.125) (5,\vara*0.0625)(6,\vara*0.03125)
%    };
%    \addplot[dashed,line width=0.75pt] coordinates { % h^2
%      (1,\vara) (2,\vara*0.5*0.5) (3,\vara*0.25*0.25) (4,\vara*0.125*0.125) (5,\vara*0.0625*0.0625)(6,\vara*0.03125*0.03125)
%    };
    \addplot[dotted,line width=0.75pt] coordinates { % h^3
    (1,\varb) (2,\varb*0.5*0.5*0.5) (3,\varb*0.25*0.25*0.25) (4,\varb*0.125*0.125*0.125) (5,\varb*0.0625*0.0625*0.0625)(6,\varb*0.03125*0.03125*0.03125)
    };
     \addplot[dashed,line width=0.75pt] coordinates { % h^4
    (1,\vara) (2,\vara*0.5*0.5*0.5*0.5) (3,\vara*0.25*0.25*0.25*0.25) (4,\vara*0.125*0.125*0.125*0.125) (5,\vara*0.0625*0.0625*0.0625*0.0625)(6,\vara*0.03125*0.03125*0.03125*0.03125)
    };
    \legend{$e_{L^2}^{\bu}$, $e_{PL^2}^{\bu}$, $e_{H^1}^{\bu}$, $e_{A}^{\bu}$, $e_{M}^{p}$, $\mathcal{O}(h^{3})$, $\mathcal{O}(h^{4})$}
  \end{semilogyaxis}
  \end{tikzpicture}
  }
%    \caption{Error norms for the sphere and $k=l=1$ with $\rho=\tilde\rho=1$.}
% \end{figure}       
\end{minipage}
\caption{Consistent formulation \eqref{discreteform1} on the torus} \label{tablePh1torus}
\end{figure}

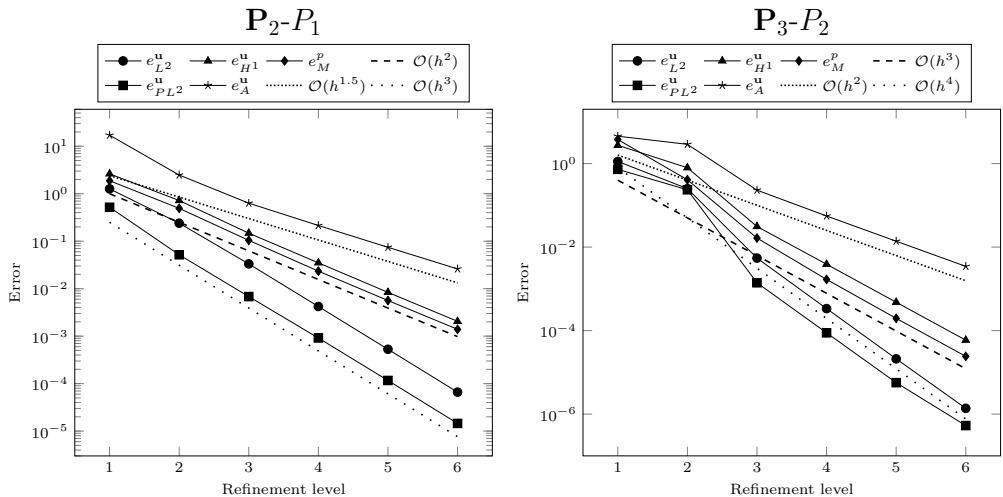
\begin{figure}
\hspace*{-0.7cm}
\begin{minipage}{0.51\textwidth}
\scriptsize
      % \begin{figure}  
\scalebox{0.81}{
  \begin{tikzpicture}
  \def\vara{1.0}
  \def\varb{0.25}
  \def\varc{2.4}
\begin{semilogyaxis}[ xlabel={Refinement level}, ylabel={Error}, ymin=3E-6, ymax=60, legend style={ at={(0.5,1.02)}, anchor=south, legend columns =2, transpose legend }, legend cell align=left, cycle list name=mark list, title= {$\vect P_2$-$P_{1}$}, title style={at={(0.5,1.16)}, align=center, font=\Large} ]
    \addplot table[x=level, y=L2] {Ph2P2P1torusopt.dat};
    \addplot table[x=level, y=L2P] {Ph2P2P1torusopt.dat};
    \addplot table[x=level, y=realH1semi] {Ph2P2P1torusopt.dat};
    \addplot table[x=level, y=uerror] {Ph2P2P1torusopt.dat};
    \addplot table[x=level, y=merr] {Ph2P2P1torusopt.dat};
%    \addplot[dashed,line width=0.75pt] coordinates { % h
%    (1,\vara) (2,\vara*0.5) (3,\vara*0.25) (4,\vara*0.125) (5,\vara*0.0625)(6,\vara*0.03125)
%    };
    \addplot[densely dotted,line width=0.75pt] coordinates { % h^1.5
      (1,\varc) (2,\varc*0.5*0.707106781) (3,\varc*0.25*0.5) (4,\varc*0.125*0.353553391) (5,\varc*0.0625*0.25)(6,\varc*0.03125*0.176776695)
    };
    \addplot[dashed,line width=0.75pt] coordinates { % h^2
      (1,\vara) (2,\vara*0.5*0.5) (3,\vara*0.25*0.25) (4,\vara*0.125*0.125) (5,\vara*0.0625*0.0625)(6,\vara*0.03125*0.03125)
    };
    \addplot[loosely dotted,line width=0.75pt] coordinates { % h^3
    (1,\varb) (2,\varb*0.5*0.5*0.5) (3,\varb*0.25*0.25*0.25) (4,\varb*0.125*0.125*0.125) (5,\varb*0.0625*0.0625*0.0625)(6,\varb*0.03125*0.03125*0.03125)
    };
    \legend{$e_{L^2}^{\bu}$, $e_{PL^2}^{\bu}$, $e_{H^1}^{\bu}$, $e_{A}^{\bu}$, $e_{M}^{p}$, $\mathcal{O}(h^{1.5})$, $\mathcal{O}(h^{2})$, $\mathcal{O}(h^{3})$}
  \end{semilogyaxis}
  \end{tikzpicture}
  }
%    \caption{Error norms for the sphere and $k=l=1$ with $\rho=\tilde\rho=1$.}
% \end{figure}    
\end{minipage}
\begin{minipage}{0.51\textwidth}
\scriptsize
      % \begin{figure}  
\scalebox{0.81}{
  \begin{tikzpicture}
  \def\vara{0.8}
  \def\varb{0.4}
  \def\varc{1.6}
\begin{semilogyaxis}[ xlabel={Refinement level}, ylabel={Error}, ymin=1E-7, ymax=20, legend style={ at={(0.5,1.02)}, anchor=south, legend columns =2, transpose legend }, legend cell align=left, cycle list name=mark list, title= {$\vect P_3$-$P_{2}$}, title style={at={(0.5,1.16)}, align=center, font=\Large} ]
    \addplot table[x=level, y=L2] {Ph2P3P2torusopt.dat};
    \addplot table[x=level, y=L2P] {Ph2P3P2torusopt.dat};
    \addplot table[x=level, y=realH1semi] {Ph2P3P2torusopt.dat};
    \addplot table[x=level, y=uerror] {Ph2P3P2torusopt.dat};
    \addplot table[x=level, y=merr] {Ph2P3P2torusopt.dat};
%    \addplot[dashed,line width=0.75pt] coordinates { % h
%    (1,\vara) (2,\vara*0.5) (3,\vara*0.25) (4,\vara*0.125) (5,\vara*0.0625)(6,\vara*0.03125)
%    };
    \addplot[densely dotted,line width=0.75pt] coordinates { % h^2
      (1,\varc) (2,\varc*0.5*0.5) (3,\varc*0.25*0.25) (4,\varc*0.125*0.125) (5,\varc*0.0625*0.0625)(6,\varc*0.03125*0.03125)
    };
    \addplot[dashed,line width=0.75pt] coordinates { % h^3
    (1,\varb) (2,\varb*0.5*0.5*0.5) (3,\varb*0.25*0.25*0.25) (4,\varb*0.125*0.125*0.125) (5,\varb*0.0625*0.0625*0.0625)(6,\varb*0.03125*0.03125*0.03125)
    };
     \addplot[loosely dotted,line width=0.75pt] coordinates { % h^4
    (1,\vara) (2,\vara*0.5*0.5*0.5*0.5) (3,\vara*0.25*0.25*0.25*0.25) (4,\vara*0.125*0.125*0.125*0.125) (5,\vara*0.0625*0.0625*0.0625*0.0625)(6,\vara*0.03125*0.03125*0.03125*0.03125)
    };
    \legend{$e_{L^2}^{\bu}$, $e_{PL^2}^{\bu}$, $e_{H^1}^{\bu}$, $e_{A}^{\bu}$, $e_{M}^{p}$, $\mathcal{O}(h^{2})$, $\mathcal{O}(h^{3})$, $\mathcal{O}(h^{4})$}
  \end{semilogyaxis}
  \end{tikzpicture}
  }
%    \caption{Error norms for the sphere and $k=l=1$ with $\rho=\tilde\rho=1$.}
% \end{figure}       
\end{minipage}
\caption{Inconsistent formulation \eqref{discreteform2} on the torus}\label{tablePh2torus}
\end{figure}

For the torus we performed  experiments for the consistent and inconsistent method and with $\vect P_2$-$P_1$ and $\vect P_3$-$P_2$ finite elements. Again we used the default parameters. In Figure \ref{tablePh1torus} we show the results for the consistent formulation \eqref{discreteform1} and in Figure \ref{tablePh2torus}  for the inconsistent formulation \eqref{discreteform2}. The observed convergence rates are  the same as for the sphere. For $\vect P_3$-$P_{2}$ finite elements we observe for both formulations a slight deterioration of the convergence rate for the tangential $L^2$-error $e_{PL^2}^{\bu}$ in the last refinement step. This  may be due to a relatively  large condition number of the stiffness matrix.

\section{Conclusions and outlook} \label{sectionconclusion}

We proposed two trace finite element methods for discretization of the surface Stokes equation. Both methods use the same penalty approach for treating the tangential flow constraint and the same Taylor-Hood $\vect P_k$--$P_{k-1}$ spaces. For a higher order geometry approximation the parametric trace finite element technique is  used. For the parameters in these methods specific choices are proposed.   The numerical experiments show that for $k=2$ and $k=3$ the resulting  methods have optimal convergence orders in the $H^1(\Gamma_h)$- and $L^2(\Gamma_h)$-norm.   For the consistent method \eqref{discreteform1}  an approximation of the Weingarten map has to be determined, which is not needed in the inconsistent method \eqref{discreteform2}.  For the consistent method and $k=2$ an optimal order discretization error bound for the case $\Gamma_h=\Gamma$ (i.e., no geometry errors) is derived in  \cite{OlshanskiiZhiliakov2019}. For the inconsistent method a rigorous optimal error bound is not available, yet. 
 
In future work these methods will be applied to other related problem classes, e.g., time dependent surface Navier-Stokes equations, and compared to other methods. Furthermore, the  analysis can  be extended in several directions, for example,  by including   geometry errors and deriving (optimal) error bounds also for the inconsistent method.

%\newpremcomparison}age
%Summary: In the experiments we saw that applying the results of the analysis of the surface vector-Laplace problem (\cite{jankuhn2019}) and stabilizing the pressure with the normal gradient stabilization (choosing $\tilde{\rho} = h$) leads to a method with optimal convergence orders in all the examined errors for the consistent formulation \eqref{discreteform1} and all but the energy norm $e_{A_{h}}^{\bu}$ for the inconsistent formulation \eqref{discreteform2}. For the special case of the sphere and a divergence-free solution, i.e. $g=0$, the inconsistent formulation \eqref{discreteform2} results in optimal convergence orders for the energy norm $e_{A_{h}}^{\bu}$ as well.

\bibliographystyle{siam}
\bibliography{literatur}{}

\begin{thebibliography}{10}

\bibitem{ngsolve}
{\em {Netgen/NGSolve}}.
\newblock https://ngsolve.org/.

\bibitem{arroyo2009}
{\sc M.~Arroyo and A.~DeSimone}, {\em Relaxation dynamics of fluid membranes},
  Phys. Rev. E, 79 (2009), p.~031915.

\bibitem{Bonito2019a}
{\sc A.~Bonito, A.~Demlow, and M.~Licht}, {\em A divergence-conforming finite
  element method for the surface {Stokes} equation}, arXiv:1908.11460,  (2019).

\bibitem{Bonito2019}
{\sc A.~Bonito, A.~Demlow, and R.~H. Nochetto}, {\em Finite element methods for
  the {Laplace-Beltrami} operator}, arXiv:1906.02786,  (2019).

\bibitem{cutFEM}
{\sc E.~Burman, S.~Claus, P.~Hansbo, M.~G. Larson, and A.~Massing}, {\em
  Cutfem: Discretizing geometry and partial differential equations},
  International Journal for Numerical Methods in Engineering, 104 (2015),
  pp.~472--501.

\bibitem{burmanembedded}
{\sc E.~Burman, P.~Hansbo, M.~G. Larson, and A.~Massing}, {\em Cut finite
  element methods for partial differential equations on embedded manifolds of
  arbitrary codimensions}, ESAIM: Mathematical Modelling and Numerical
  Analysis, 52 (2018), pp.~2247--2282.

\bibitem{demlow2009higher}
{\sc A.~Demlow}, {\em Higher-order finite element methods and pointwise error
  estimates for elliptic problems on surfaces}, SIAM J. Numer. Anal., 47
  (2009), pp.~805--827.

\bibitem{DEreview}
{\sc G.~Dziuk and C.~M. Elliott}, {\em Finite element methods for surface
  {PDEs}}, Acta Numerica, 22 (2013), pp.~289--396.

\bibitem{ebin1970groups}
{\sc D.~G. Ebin and J.~Marsden}, {\em Groups of diffeomorphisms and the motion
  of an incompressible fluid}, Annals of Mathematics, 92 (1970), pp.~102--163.

\bibitem{fries2018higher}
{\sc T.-P. Fries}, {\em Higher-order surface {FEM} for incompressible
  {Navier-Stokes} flows on manifolds}, International Journal for Numerical
  Methods in Fluids, 88 (2018), pp.~55--78.

\bibitem{grande2017higher}
{\sc J.~Grande, C.~Lehrenfeld, and A.~Reusken}, {\em Analysis of a high-order
  trace finite element method for {PDEs} on level set surfaces}, SIAM Journal
  on Numerical Analysis, 56 (2018), pp.~228--255.

\bibitem{GurtinMurdoch75}
{\sc M.~E. Gurtin and A.~I. Murdoch}, {\em A continuum theory of elastic
  material surfaces}, Archive for Rational Mechanics and Analysis, 57 (1975),
  pp.~291--323.

\bibitem{hansbo2016analysis}
{\sc P.~Hansbo, M.~G. Larson, and K.~Larsson}, {\em Analysis of finite element
  methods for vector {Laplacians} on surfaces}, IMA Journal of Numerical
  Analysis,  (2019).

\bibitem{Jankuhn1}
{\sc T.~Jankuhn, M.~A. Olshanskii, and A.~Reusken}, {\em Incompressible fluid
  problems on embedded surfaces: Modeling and variational formulations},
  Interfaces and Free Boundaries, 20 (2018), pp.~353--377.

\bibitem{jankuhn2019}
{\sc T.~Jankuhn and A.~Reusken}, {\em Trace finite element methods for surface
  vector-{Laplace} equations}, arXiv:1904.12494,  (2019).

\bibitem{Kobaetal_QAM_2017}
{\sc H.~Koba, C.~Liu, and Y.~Giga}, {\em {Energetic variational approaches for
  incompressible fluid systems on an evolving surface}}, Quart. Appl. Math., 75
  (2017), pp.~359--389.

\bibitem{Lederer2019}
{\sc P.~L. Lederer, C.~Lehrenfeld, and J.~Sch\"oberl}, {\em Divergence-free
  tangential finite element methods for incompressible flows on surfaces},
  arXiv:1909.06229,  (2019).

\bibitem{ngsxfem}
{\sc C.~Lehrenfeld}, {\em ngsxfem}.
\newblock https://github.com/ngsxfem.

\bibitem{lehrenfeld2016high}
{\sc C.~Lehrenfeld}, {\em High order unfitted finite element methods on level
  set domains using isoparametric mappings}, Computer Methods in Applied
  Mechanics and Engineering, 300 (2016), pp.~716--733.

\bibitem{lehrenfeld2017analysis}
{\sc C.~Lehrenfeld and A.~Reusken}, {\em Analysis of a high-order unfitted
  finite element method for elliptic interface problems}, IMA Journal of
  Numerical Analysis, 38 (2017), pp.~1351--1387.

\bibitem{mitrea2001navier}
{\sc M.~Mitrea and M.~Taylor}, {\em {N}avier-{S}tokes equations on {Lipschitz}
  domains in {Riemannian} manifolds}, Mathematische Annalen, 321 (2001),
  pp.~955--987.

\bibitem{miura2017singular}
{\sc T.-H. Miura}, {\em On singular limit equations for incompressible fluids
  in moving thin domains}, Quart. Appl. Math., 76 (2018), pp.~215--251.

\bibitem{Nitschkeetal_arXiv_2018}
{\sc I.~Nitschke, S.~Reuther, and A.~Voigt}, {\em Hydrodynamic interactions in
  polar liquid crystals on evolving surfaces}, Phys. Rev. Fluids, 4 (2019),
  p.~044002.

\bibitem{nitschke2012finite}
{\sc I.~Nitschke, A.~Voigt, and J.~Wensch}, {\em A finite element approach to
  incompressible two-phase flow on manifolds}, Journal of Fluid Mechanics, 708
  (2012), pp.~418--438.

\bibitem{olshanskii2018finite}
{\sc M.~A. Olshanskii, A.~Quaini, A.~Reusken, and V.~Yushutin}, {\em A finite
  element method for the surface {Stokes} problem}, SIAM Journal on Scientific
  Computing, 40 (2018), pp.~A2492--A2518.

\bibitem{olshanskii2016trace}
{\sc M.~A. Olshanskii and A.~Reusken}, {\em Trace finite element methods for
  {PDEs} on surfaces}, in Geometrically Unfitted Finite Element Methods and
  Applications, S.~P.~A. Bordas, E.~Burman, M.~G. Larson, and M.~A. Olshanskii,
  eds., Cham, 2017, Springer International Publishing, pp.~211--258.

\bibitem{OlshanskiiZhiliakov2019}
{\sc M.~A. Olshanskii, A.~Reusken, and A.~Zhiliakov}, {\em Inf-sup stability of
  the trace {$P_2$-$P_1$ {Taylor-Hood} elements for surface {PDEs}}}, Preprint
  arXiv:1909.02990,  (2019).

\bibitem{olshanskii2019penalty}
{\sc M.~A. Olshanskii and V.~Yushutin}, {\em A penalty finite element method
  for a fluid system posed on embedded surface}, Journal of Mathematical Fluid
  Mechanics, 21 (2019), p.~14.

\bibitem{reusken2015analysis}
{\sc A.~Reusken}, {\em Analysis of trace finite element methods for surface
  partial differential equations}, IMA Journal of Numerical Analysis, 35
  (2015), pp.~1568--1590.

\bibitem{reusken2018stream}
\leavevmode\vrule height 2pt depth -1.6pt width 23pt, {\em Stream function
  formulation of surface {Stokes} equations}, IMA Journal of Numerical
  Analysis,  (2018).

\bibitem{reuther2015interplay}
{\sc S.~Reuther and A.~Voigt}, {\em The interplay of curvature and vortices in
  flow on curved surfaces}, Multiscale Modeling \& Simulation, 13 (2015),
  pp.~632--643.

\bibitem{reuther2018solving}
\leavevmode\vrule height 2pt depth -1.6pt width 23pt, {\em Solving the
  incompressible surface {Navier-Stokes} equation by surface finite elements},
  Physics of Fluids, 30 (2018), p.~012107.

\bibitem{Schoeberl1997}
{\sc J.~Sch{\"o}berl}, {\em Netgen an advancing front 2d/3d-mesh generator
  based on abstract rules}, Computing and Visualization in Science, 1 (1997),
  pp.~41--52.

\bibitem{taylor1992analysis}
{\sc M.~E. Taylor}, {\em Analysis on {Morrey} spaces and applications to
  {Navier-Stokes} and other evolution equations}, Communications in Partial
  Differential Equations, 17 (1992), pp.~1407--1456.

\bibitem{Temam88}
{\sc R.~Temam}, {\em Infinite-dimensional dynamical systems in mechanics and
  physics}, Springer, New York, 1988.

\end{thebibliography}

\end{document}